\documentclass{article}
\usepackage[utf8]{inputenc}
\usepackage{amsfonts, amsmath, amssymb, amsthm, authblk, graphicx, enumitem}
\usepackage{cite}
\usepackage[colorlinks=true,urlcolor=blue, citecolor=red,linkcolor=blue]{hyperref}
\usepackage[left=2cm, right=2cm, top=3cm]{geometry}
\usepackage[capitalise]{cleveref}

\newtheoremstyle{sltheorem}
{}                
{}                
{\itshape}        
{}                
{\bfseries}       
{.}               
{ }               
{}				  
\theoremstyle{sltheorem}
\newtheorem{teorema}{Theorem}[section]

\newtheorem{lema}[teorema]{Lemma}
\newtheorem{proposicion}[teorema]{Proposition}

\theoremstyle{definition}

\newtheorem{observacion}[teorema]{Remark}

\numberwithin{equation}{section}

\renewcommand{\leq}{\leqslant}

\renewcommand{\geq}{\geqslant}

\newcommand*\diff{\mathop{}\!\mathrm{d}}

\begin{document}
\title{Controllability results for cascade systems of $m$ coupled $N$-dimensional Stokes and Navier-Stokes systems by $N-1$ scalar controls}
\author[1]{Tak\'eo Takahashi \thanks{takeo.takahashi@inria.fr }}
\affil[1]{Universit\'e de Lorraine, CNRS, Inria, IECL, F-54000 Nancy, France}
\author[2]{Luz de Teresa \thanks{ldeteresa@im.unam.mx. Partially supported by Conacyt, project A1-S-17475}}
\affil[2]{Instituto de Matematicas, Universidad Nacional Aut\'onoma de M\'exico, Circuito Exterior, C.U., C. P. 04510 Ciudad de México, México}
\author[2]{Ying Wu-Zhang\thanks{yingwu@im.unam.mx}}

\maketitle
	
\tableofcontents

\begin{abstract}
In this paper we deal with the controllability properties of a system of $m$ coupled Stokes systems or $m$ coupled Navier-Stokes systems. 
We show the null-controllability of such systems in the case where the coupling is in a cascade form and when the control acts only on one of the systems. Moreover, we impose that this control has a vanishing component so that we control a $m\times N$ state (corresponding to the velocities of the fluids) by  $N-1$ distributed scalar controls. The proof of the controllability of the coupled Stokes system is based on a Carleman estimate for the adjoint system. The local null-controllability of the coupled Navier-Stokes systems is then obtained by means of the source term method and a Banach fixed point.
\end{abstract}

\vspace{1cm}



\section{Introduction}

Controllability  issues related to a single parabolic equation or to a single Stokes or Navier-Stokes system have been intensively studied in the last fifty years giving rise to interesting techniques, new challenges and open problems. See some seminal results \cite{fursikov1996controllability, lebeau1995, fattorini1975} for the heat equation and \cite{fernandez2004local, MR1470445, imanuvilov1998} for the Navier-Stokes system.
The literature is vast and it is difficult to mention all the intensive studies about this subject. 
However, it is only in the last fifteen years that the challenging issue
of controlling coupled parabolic systems has attracted the interest of the control
community. This kind of systems appears mathematically  in optimal control theory as a characterization of the optimal control (with one equation coupled to its adjoint) but also appears, for example, in the study of chemical
reactions (see e.g. \cite{ErdiToth1989,Conforto2017}), and in a wide variety of mathematical biology and
physical situations (see e.g.\cite{ Iida2017} ). In the case of scalar (heat) coupled equations an important number of challenging problems has been solved (see \cite{Ammarkhodja2011} for a survey of results until 2011) and sometimes the results have been surprising \cite{Ammarkhodja2014,Ammarkhodja2016, Benabdallah2020}. In the case of coupled Stokes or Navier-Stokes systems, to our knowledge, only some cases of \emph{two} coupled systems have been treated \cite{Guerrero2007, Carreno-Gueye2014, MontoyadeTeresa2018, MR3348416}. 
Here our aim is to generalize results for a $m$ scalar cascade system \cite{gonzalez2010controllability}  to a $m$ $N$-dimensional Stokes or Navier-Stokes cascade system but including an extra deal: to eliminate one component on the $N$-dimensional control.

Let us be more specific: we consider a bounded domain $\Omega$  of $\mathbb{R}^N$ ($N=2, 3$) whose boundary $\partial\Omega$ is regular enough. 
Let $T>0$ and let $\omega\subset\Omega$ be a (arbitrary small) nonempty open subset which will usually be referred as the \textit{control domain}. 
We will use the notation $Q=\Omega\times(0,T)$ and $\Sigma=\partial\Omega\times(0,T).$ 
	
In this article, we are interested in the null controllability of a coupled system of $m$ Stokes or Navier-Stokes systems, with $m\geq 2$:
$$
\left\{\begin{array}{lll}
	\partial_t y^{(i)}-\nu_i\Delta y^{(i)}+\varepsilon \left(y^{(i)}\cdot \nabla\right)y^{(i)}+\nabla p^{(i)} = \sum_{j=1}^m \left(B_{i,j}\cdot\nabla\right) y^{(j)}+\sum_{j=1}^m A_{i,j} y^{(j)}+D_{i} v 1_{\omega}  &\text{in} \ Q, &(1\leq i\leq m)\\
	\nabla\cdot y^{(i)}=0 &\text{in }Q, &(1\leq i\leq m)\\
	y^{(i)}=0 &\text{on }\Sigma, &(1\leq i\leq m)\\
	y^{(i)}(\cdot,0)=y_0^{(i)} &\text{in }\Omega, &(1\leq i\leq m)
\end{array}\right.
$$
with $\varepsilon=0$ for Stokes systems and $\varepsilon=1$ for Navier-Stokes systems, where $A_{i,j}\in \mathcal{M}_N(\mathbb{R})$, $B_{i,j}\in \mathbb{R}^N$ and 
where $D_i\in \mathcal{M}_{N,r}(\mathbb{R})$ for some $r\in \mathbb{N}^*$.
We have denoted by $1_\omega$ the characteristic function of $\omega$. 
The constants $\nu_i>0$ are the viscosities of the fluids.

We can write the above systems in a more compact way as
\begin{equation} \label{eq:sistemamatriz}
	\begin{cases}
		\partial_t y-\nu\Delta y+\nabla p=\left(B\cdot \nabla\right)y+Ay+Dv1_{\omega} &\text{in } Q, \\
		\nabla\cdot y=0 &\text{in } Q, \\
		y=0 &\text{on }\Sigma, \\
		y\left(\cdot,0\right)=y_0 &\text{in }\Omega,
	\end{cases}
\end{equation}
or
\begin{equation} \label{eq:sistemamatriz+navier}
	\begin{cases}
		\partial_t y-\nu\Delta y+(y\cdot\nabla)y+\nabla p=\left(B\cdot \nabla\right)y+ Ay+Dv1_{\omega} &\text{in } Q, \\
		\nabla\cdot y=0 &\text{in } Q, \\
		y=0 &\text{on }\Sigma, \\
		y\left(\cdot,0\right)=y_0 &\text{in }\Omega,
	\end{cases}
\end{equation}
where we set
$$
y=\begin{pmatrix}
	y^{(i)}
\end{pmatrix}_{1\leq i\leq m},
\quad
\nu\Delta y=\begin{pmatrix}
	\nu_i\Delta y^{(i)}
\end{pmatrix}_{1\leq i\leq m},
\quad
\left(y\cdot\nabla\right)y=\begin{pmatrix}
	\left(y^{(i)}\cdot\nabla\right)y^{(i)}
\end{pmatrix}_{1\leq i\leq m},
\quad
\nabla p=\begin{pmatrix}
	\nabla p^{(i)}
\end{pmatrix}_{1\leq i\leq m},
$$
$$
\nabla\cdot y=\begin{pmatrix}
	\nabla\cdot y^{(i)}
\end{pmatrix}_{1\leq i\leq m},
\quad
y_0=\begin{pmatrix}
	y_0^{(i)}
\end{pmatrix}_{1\leq i\leq m}
$$
and
$$
(B\cdot\nabla) y=\begin{pmatrix}
	\sum_{j=1}^m (B_{i,j}\cdot\nabla) y^{(j)}
\end{pmatrix}_{1\leq i\leq m},
\quad
Ay=\begin{pmatrix}
\sum_{j=1}^m A_{i,j} y^{(j)}
\end{pmatrix}_{1\leq i\leq m},
\quad
Dv=\begin{pmatrix}
D_i v
\end{pmatrix}_{1\leq i\leq m}.
$$

In this work, we will focus in the particular case where the partitioned matrix $A$ has the form 
\begin{equation}\label{shapematA}
A=
\begin{pmatrix}
	A_{1,1} & A_{1,2} & A_{1,3} & \ldots & A_{1,m} \\
	A_{2,1} & A_{2,2} & A_{2,3} & \ldots & A_{2,m} \\
	0 & A_{3,2} & A_{3,3} & \ldots & A_{3,m} \\
	\vdots & \vdots & \ddots & \ddots & \vdots \\
	0 & 0 & \ldots & A_{m,m-1} & A_{m,m}
\end{pmatrix},
\quad 
	B=
	\begin{pmatrix}
		B_{1,1} & B_{1,2} & \ldots & B_{1,m} \\
		0 & B_{2,2}  & \ldots & B_{2,m} \\
		\vdots & \vdots  & \ddots & \vdots \\
		0 & 0 & \ldots & B_{m,m}
	\end{pmatrix},
\end{equation}
with all the blocks under the diagonal non zero for the matrix $A$ and such that all its blocks are scalar matrices.
More precisely, our hypotheses on $A$ and $B$ are
\begin{equation}\label{hypmatA}
	A_{i,j}=a_{i,j} I_N, \quad a_{i,i-1}\neq 0 \quad (2\leq i\leq m), \quad a_{i,j}=0 \quad \text{if} \ i\geq j+2,
\end{equation}
\begin{equation}\label{hypmatB}
	B_{i,j}=0 \quad \text{if} \ i\geq j+1.
\end{equation}
We also control \eqref{eq:sistemamatriz} or \eqref{eq:sistemamatriz+navier} by acting only on one of the Stokes or of the Navier-Stokes systems, for instance the first one, and with $N-1$ scalar controls on this system. Thus, without loss of generality, we assume
\begin{equation}\label{hypmatD}
r=N-1, \quad D_j=0 \quad (j\geq 2), \quad D_1=\begin{pmatrix}
1 \\ 0
\end{pmatrix} \quad (\text{if}\ N=2) \quad \text{or} \quad
D_1=\begin{pmatrix}
1 & 0
\\ 
0 & 1
\\
0 & 0
\end{pmatrix} \quad (\text{if}\  N=3).
\end{equation}
The above choice of the matrix $A$ corresponds to a particular coupling considered in the context of the null-controllability of systems of $m$ linear heat equations, see \cite{gonzalez2010controllability} and our aim is to extend this result in the case of coupled of Stokes or Navier-Stokes systems.

In order to state our main results, we recall some standard functional spaces associated with the Stokes system:
\begin{equation}\label{defH}
    H=\left\{y\in L^2(\Omega)^N:\nabla\cdot y=0\text{ in }\Omega,\ y\cdot \textbf{n}=0\text{ on }\partial\Omega\right\},
\end{equation}
\begin{equation}\label{defV}
    V=\left\{y\in H_0^1(\Omega)^N:\nabla\cdot y=0\text{ in }\Omega\right\}
\end{equation}
and
\begin{equation}\label{defcalH}
    \mathcal{H}:=H^m, \quad \mathcal{V}:=V^m.
\end{equation} 
Our main result is the following theorem.
\begin{teorema} \label{teo:1}
    Assume \eqref{hypmatA}-\eqref{hypmatD}. Then, for any $T>0$ and for any $y_0\in \mathcal{H}$, there exists a control $v\in L^2(0,T;L^2(\omega)^{N-1})$ 
    such that the corresponding solution $y=\left(y^{(1)},...,y^{(m)}\right)$ to \eqref{eq:sistemamatriz} satisfies 
    $$
    y(\cdot,T)=0\quad \text{in }\Omega.
    $$
\end{teorema}
\begin{observacion}
    As a consequence, we deduce that we can control the system \eqref{eq:sistemamatriz}  of $N\times m$ scalar equations with $N-1$ 
    scalar controls.
\end{observacion}
From \cref{teo:1} and a general method to deal with the controllability of nonlinear parabolic systems, we deduce 
the local null controllability of the system \eqref{eq:sistemamatriz+navier}:
\begin{teorema} \label{teo:2}
    Assume \eqref{hypmatA}-\eqref{hypmatD}. Then, for any $T>0$, there exists $\delta>0$ such that, for any $y_0\in \mathcal{V}$ satisfying 
    $$
    \|y_0\|_{\mathcal{V}}\leq\delta,
    $$
    there exists a control $v\in L^2(0,T;L^2(\omega)^{N-1})$ such that the corresponding solution $y=\left(y^{(1)},...,y^{(m)}\right)$ to \eqref{eq:sistemamatriz+navier} 
    satisfies
    $$
    y(\cdot,T)=0\quad \text{in }\Omega.
    $$
\end{teorema}
%
In order to prove \cref{teo:1}, we introduce the adjoint system of \eqref{eq:sistemamatriz}: 
\begin{equation}\label{eq:adjointmatriz}
	\begin{cases}
		-\partial_t \varphi-\Delta \varphi+\left(B^*\cdot\nabla\right)\varphi-A^*\varphi+\nabla \pi =0 &\text{in }Q, \\
		\nabla\cdot \varphi=0 &\text{in }Q, \\
		\varphi=0 & \text{on }\Sigma, \\
		\varphi(\cdot,T)=\varphi_T & \text{in }\Omega,
	\end{cases}
\end{equation}
or in its expanded form:
\begin{equation} \label{eq:systemadjoint2}
    \left\{\begin{array}{lll}
		-\partial_t \varphi^{(i)}- \Delta \varphi^{(i)}+\nabla \pi^{(i)} +\sum_{j=1}^i\left((B_{j,i}\cdot\nabla)\varphi^{(j)}-A_{j,i}\varphi^{(j)}\right)=A_{i+1,i}\varphi^{(i+1)} &\text{in }Q, & (1\leq i\leq m-1)\\
		\vdots \\
		-\partial_t \varphi^{(m)}- \Delta \varphi^{(m)}+\nabla \pi^{(m)} +\sum_{j=1}^m\left((B_{j,m}\cdot\nabla)\varphi^{(j)}-A_{j,m}\varphi^{(j)}\right)=0 &\text{in }Q,\\
		\nabla\cdot \varphi^{(i)}=0 &\text{in }Q,  &(1\leq i\leq m) \\
		\varphi^{(i)}=0 &\text{on }\Sigma, &(1\leq i\leq m) \\
		\varphi^{(i)}(\cdot,T)=\varphi^{(i)}_T & \text{in }\Omega, &(1\leq i\leq m) 
	\end{array}\right.
\end{equation}	
with $\varphi^{(i)}_T\in H$ $(1\leq i\leq m)$ and we also denote by $\varphi_j^{(i)}$, $j=1,\ldots,N$ the coordinates of $\varphi^{(i)}$.
Note that by setting 
$$
\varphi^{(m+1)} \equiv 0, \quad A_{m+1,i}=0 \quad (1\leq i\leq m),
$$
we can write the first above equations as
\begin{equation}\label{tt00}
    -\partial_t \varphi^{(i)}- \Delta \varphi^{(i)}+\nabla \pi^{(i)} +\sum_{j=1}^i\left((B_{j,i}\cdot\nabla)\varphi^{(j)}-A_{j,i}\varphi^{(j)}\right)=A_{i+1,i}\varphi^{(i+1)} \quad \text{in} \ Q, \quad (1\leq i\leq m).
\end{equation}
Following a standard duality argument (see, for instance, \cite[Theorem 11.2.1, p.357]{TucsnakWeiss}), \cref{teo:1} will be obtained as a consequence of
the following observability inequality: 
\begin{equation}\label{eq:observability}
    \sum_{i=1}^m \int_{\Omega}  \left|\varphi^{(i)}(x,0)\right|^2  \diff{x} \leq C(T)\sum_{j=1}^{N-1}\iint_{\omega\times(0,T)}  \left| \varphi_j^{(1)} \right|^2  \diff{x} \diff{t}.
\end{equation}
for some $C$ depending on $T$, $\Omega$ and $\omega$.

This observability inequality is a consequence of the Carleman inequality obtained in \cref{prop:0}.
Such a method based on \textit{(global) Carleman inequalities} for the controllability of parabolic equations
was introduced in \cite{fursikov1996controllability}. Such inequalities have been used by many authors to deal with Stokes of Navier-Stokes systems (for instance,
\cite{imanuvilov2001remarks} or \cite{fernandez2004local}). The case of controls with some vanishing components was considered in 
\cite{fernandez2006some}, \cite{coron2009null} and \cite{Carreno-Gueye2014}. 
An important work related to this subject is \cite{MR3279537} where the authors obtain the local null controllability of the Navier-Stokes system in dimension 3 with a control having two vanishing components. In that case, the method is based on a different linearization and on a different approach based on results of Gromov.
We follow here the method introduced in \cite{coron2009null}.
As a first step, one can get rid of the pressure by applying a differential operator on \eqref{tt00} (or on components of \eqref{tt00}) such as 
$\operatorname{curl}$ or $\Delta$. This leads to a system of coupled heat equations but without prescribed boundary conditions. Using results such as \cite{MR2224822} or \cite{imanuvilov2009carleman}, one can obtain a Carleman estimate with boundary terms that can be absorbed by some standard arguments.

Let us point out that here we consider the operator $\nabla^2 \Delta$ to get rid of the pressure. In the case of Navier slip boundary conditions, the authors of \cite{SergioNavier} and \cite{SergioCristhian} only consider the operator $\nabla \Delta$. Due to the coupling between the Stokes systems, it was more convenient to use the operator $\nabla^2 \Delta$. 
Another important remark is that to recover the observability on the components that are not observed, one has to use the divergence-free condition on $\varphi^{(i)}$ and the Dirichlet boundary conditions. Unhappily in this process, one loses part of the weights on these components, as it can be seen in the definition of the weights in \cref{prop:0}
(see \eqref{defI1}). 

The last part of the proof of \eqref{eq:observability}  consists in estimating the local terms associated with $\varphi_j^{(i)}$, $i>1$ and this is done by using \eqref{hypmatA}.
It is important to notice that in the proof of the Carleman estimate of \cref{prop:0}, we consider the case where the adjoint system \eqref{eq:adjointmatriz} has no right-hand side. Due to the method of proof, with the use of the operator $\nabla^2 \Delta$, this would impose restrictions on the regularity of the source term, see for instance \cite{Carreno-Gueye2014} where the authors consider the coupling between two Stokes systems and where one of the source needs to be $H^1$ in space. Nevertheless, using the general method introduced in \cite{liu2013single}, we can use \eqref{eq:observability} for the adjoint system without source terms to deal with the controllability of the Navier-Stokes systems \eqref{eq:sistemamatriz+navier}.


This paper is organized as follows. In \cref{S2}, we introduce the weights for our Carleman estimates and we recall some results, in particular some Carleman estimates for other systems. \cref{S3} corresponds to the statement and to the proof of the Carleman estimate for \eqref{eq:systemadjoint2}. 
Finally, in \cref{S4}, we use this estimate to prove \cref{teo:1} and \cref{teo:2}.

\section{Preliminaries}\label{S2}

In the whole article, we use the notation $C$ for a generic positive constant that depends on $\Omega,\omega$. 
We can assume to simplify that $T\in (0,1)$ and that will allow us to avoid some dependence on $T$ for some constants. 
We also assume that $\nu_i=1$ for $i=1,\ldots,m$ since these constants do not play any role in the analysis.
Finally we only write the proof in the case
$N=2$, the case $N=3$ can be done similarly.

\subsection{Definition of the weights and first Carleman estimates}
To write our Carleman inequalities, we introduce standard weights and functions. First, we consider 
a nonempty domain $\omega_0$ such that $\overline{\omega_0}\subset\omega$. Then, using \cite{fursikov1996controllability} 
(see also \cite[Theorem 9.4.3, p.299]{TucsnakWeiss}), there exists
$\eta^0\in C^2(\overline{\Omega})$ satisfying
	\begin{equation*}
		\eta^0>0 \ \text{in}\ \Omega, \quad 
		\eta^0= 0 \ \text{on} \ \partial\Omega,\quad
		\max_\Omega \eta^0=1,\quad
		\nabla\eta^0\neq 0 \ \text{in}\ \overline{\Omega\setminus \omega_0}.
	\end{equation*}
Then, we define the following functions:
	\begin{equation} \label{eq:weights}
		\alpha(x,t)=\frac{\exp\left\{\lambda(2\ell+2)\right\}-\exp\{\lambda (2\ell+\eta^0(x))\}}{t^\ell(T-t)^\ell},
		\quad
		\xi(x,t)=\frac{\exp\{\lambda (2\ell+\eta^0(x))\}}{t^\ell(T-t)^\ell},
	\end{equation}
	\begin{equation}\label{eq:weights2}
		\alpha^*(t)=\max_{x\in\overline{\Omega}}\alpha(x,t)=\frac{\exp\left\{\lambda(2\ell+2)\right\}-\exp\{ 2\lambda\ell\}}{t^\ell(T-t)^\ell},
		\quad 		
		\xi^*(t)=\min_{x\in\overline{\Omega}}\xi(x,t)=\frac{\exp\{2\lambda \ell\}}{t^\ell(T-t)^\ell},
	\end{equation}
where $\ell\geq 14$, $\lambda>1$. In literature $\ell>3$ is enough. However, to prove our result we need to use $\ell\geq 14$ (see Proof of Proposition \ref{prop:1}).
	
Note that we have the following useful relations: there exists $C>0$ depending on $\Omega$ such that
	\begin{equation}\label{09:48}
		\left|\partial_t\alpha\right|+\left|\partial_t \xi\right|
		\leq C T \xi^{1+1/\ell},
	\end{equation}
	\begin{equation}\label{09:49}
		\left|\left(\alpha^*\right)'\right|+\left|\left(\xi^*\right)'\right|
		\leq C T \left(\xi^*\right)^{1+1/\ell},
		\quad
		\left|\left(\alpha^*\right)''\right|+\left|\left(\xi^*\right)''\right|
		\leq C T^2 \left(\xi^*\right)^{1+2/\ell},
	\end{equation}
	\begin{equation}\label{09:50}
		\xi^*\geq \frac{C}{T^{2\ell}},
	\end{equation}
	\begin{equation}\label{11:01}
		\left|\nabla \alpha\right|=\left|\nabla \xi\right|\leq C\lambda \xi,
		\quad
		\left|\Delta \alpha\right|=\left|\Delta \xi\right|\leq C\lambda^2 \xi.
	\end{equation}

Weights of the kind \eqref{eq:weights}  were first considered in \cite{fursikov1996controllability}. In its present form, these weights have already been used in \cite{guerrero2007null} in order to obtain Carleman estimates for the controllability of strongly coupled parabolic equations and later in \cite{Guerrero2007} for the existence of insensitizing controls for Stokes systems.
	
Now, we recall some standard results. The first one is a Carleman estimate for the gradient operator, it is stated and proved in \cite{coron2009null}:
	\begin{lema} \label{lemma:21}
		Let $r\in\mathbb{R}.$ There exists $C>0$ depending only on $r,\Omega$ and $\omega_0$ such that, for every $T>0$, $s\geq C$ and every $u\in L^2(0,T;H^1(\Omega)),$
		$$
		\iint_Qe^{-2s\alpha}s^{r+2}\xi^{r+2} |u|^2 \diff{x}\diff{t}
		\leq C
		\left(
		\iint_Qe^{-2s\alpha}s^r \xi^r|\nabla u|^2 \diff{x}\diff{t}+\iint_{\omega_0\times(0,T)}e^{-2s\alpha}s^{r+2}\xi^{r+2}|u|^2 \diff{x}\diff{t}.
		\right).
		$$	
	\end{lema}
The second result is a Carleman estimate for the Laplace operator, it is stated and proved in \cite{coron2009null}:
	\begin{lema} \label{lemma:215}
		Let $r\in\mathbb{R}.$ There exists $C>0$ depending only on $r,\Omega$ and $\omega_0$ such that, for every $T>0$, $s\geq C$ and every $u\in L^2(0,T;H^2(\Omega)\cap H^1_0(\Omega)),$
\begin{multline*}
		\iint_Qe^{-2s\alpha}s^{r+3}\xi^{r+3} |u|^2 \diff{x}\diff{t}
		+\iint_Qe^{-2s\alpha}s^{r+1}\xi^{r+1} |\nabla u|^2 \diff{x}\diff{t}
		\\
		\leq C
		\left(
		\iint_Qe^{-2s\alpha}s^r \xi^r|\Delta u|^2 \diff{x}\diff{t}+\iint_{\omega_0\times(0,T)}e^{-2s\alpha}s^{r+3}\xi^{r+3}|u|^2 \diff{x}\diff{t}.
		\right).
\end{multline*}
	\end{lema}
The third result is a Carleman estimate for the heat equation with non-homogeneous Dirichlet boundary conditions. It is proved in \cite{imanuvilov2009carleman}.
	\begin{lema} \label{lemma:22}
		There exists a constant $C>0$ such that for any $\lambda \geq C$, $s\geq C$, $f_0, f_1,...,f_N \in L^2(Q)$ and 
		$$
		u\in L^2\left(0,T;H^1(\Omega)\right)\cap H^1\left(0,T;H^{-1}(\Omega)\right)
		$$ 
		satisfying 
		$$
		\partial_t u-\Delta u=f_0+\sum_{j=1}^N\partial_j f_j \quad\text{in }Q,
		$$
		we have
		\begin{multline} \label{eq:parabolic}
			\iint_Q  e^{-2s\alpha}\big(s^{-1}\xi^{-1}\left|\nabla u\right|^2+s\xi\left|u\right|^2\big)\diff{x}\diff{t}  
			\leq C\left(\iint_{\omega_0\times(0,T)}e^{-2s\alpha}s\xi \left|u\right|^2 \diff{x}\diff{t}
			\right. \\ \left.
			+\left\|e^{-s\alpha^*}s^{-1/4}(\xi^*)^{-1/4+1/\ell}u\right\|^2_{L^2(\Sigma)}
			+\left\|e^{-s\alpha^*}s^{-1/4}(\xi^*)^{-1/4}u\right\|^{2}_{H^{\frac{1}{4},\frac{1}{2}}(\Sigma)} 
			\right. \\ \left.
			+\iint_Q e^{-2s\alpha}s^{-2}\xi^{-2} \left|f_0\right|^2 \diff{x}\diff{t}+\sum_{j=1}^N \iint_Q e^{-2s\alpha}\left|f_j\right|^2\diff{x}\diff{t} \right).
		\end{multline}
	\end{lema}
Recall that
$$
\left\|u\right\|_{H^{\frac{1}{4},\frac{1}{2}}(\Sigma)}=\left(\left\|u\right\|_{H^{1/4}(0,T;L^2(\partial\Omega))}^2+\left\|u\right\|^2_{L^2(0,T;H^{1/2}(\partial\Omega))}\right)^{1/2}.
$$

\subsection{Regularity results for the Stokes systems}\label{sec_reg}
We recall that $H$, $V$, $\mathcal{H}$ and $\mathcal{V}$ are defined by \eqref{defH}, \eqref{defV} and \eqref{defcalH}. 
We denote by $P_0 : L^2(\Omega)^N \to H$ the Leray projector and we consider the projection $\mathcal{P}$ defined by:
$$
\mathcal{P} : \left[L^2(\Omega)^N\right]^m \to \mathcal{H}, \quad y=\begin{pmatrix}
y^{(1)},\ldots,y^{(m)}
\end{pmatrix}
    \mapsto 
\begin{pmatrix}
P_0y^{(1)},\ldots,P_0 y^{(m)}
\end{pmatrix}.
$$
We also consider the unbounded operators in $\mathcal{H}$
defined by
\begin{equation}\label{defcalAdom}
        \mathcal{D}(\mathcal{A}_0):=\mathcal{D}(\mathcal{A}_1):=\mathcal{D}(\mathcal{A}):=\left[H^2(\Omega)^N\cap V\right]^m,
\end{equation}
\begin{equation}\label{defcalA}
\mathcal{A}_0y:= -\mathcal{P} \Delta y, \quad \mathcal{A}_1 y:=-\mathcal{P}\left[\left(B\cdot \nabla\right)y+ Ay\right], 
\quad
\mathcal{A}:= \mathcal{A}_0+\mathcal{A}_1.
\end{equation}
It is well-known (see, for instance, \cite[Theorem 2.1.1, pp.128--129]{Sohr}) that $\mathcal{A}_0$ is self-adjoint and positive and
one can check that
\begin{equation}\label{defcalAstar}
\mathcal{D}(\mathcal{A}^*)=\left[H^2(\Omega)^N\cap V\right]^m,
\quad
\mathcal{A}^*= \mathcal{A}_0 + \mathcal{A}_1^*,
\end{equation}
where
$$
\mathcal{A}_1^* \varphi:=\mathcal{P}\left[\left(B^*\cdot \nabla\right)\varphi- A^* \varphi\right],
$$
with the notation
$$
\left(A^*\right)_{i,j}:=A_{j,i}^{\top}\in \mathcal{M}_N(\mathbb{R}), \quad \left(B^*\right)_{i,j}:=B_{j,i}\in \mathbb{R}^N.
$$
Here we have used the notation $\cdot^\top$ for the transpose of a matrix in $\mathcal{M}_N(\mathbb{R})$.
One can check that
$$
\left\| \mathcal{A}_0^{1/2} \varphi \right\|_{\mathcal{H}}=\left\| \nabla \varphi \right\|_{L^2(\Omega)^{N^2m}},
\quad
\mathcal{V}=\mathcal{D}\left(\mathcal{A}_0^{1/2}\right), 
\quad
\quad \left\| \mathcal{A}_1^* \varphi \right\|_{\mathcal{H}} \leq C \left\| \mathcal{A}_0^{1/2} \varphi \right\|_{\mathcal{H}}
$$
for some constant $C$. In particular, one deduce (see, for instance, \cite[Theorem 2.1, p.80]{Pazy}) that $-\mathcal{A}$ is the infinitesimal generator of an analytic semigroup. 
Let us now consider the adjoint system:
\begin{equation}\label{Stokes}
	\begin{cases}
		-\partial_t \varphi-\Delta \varphi+\left(B^*\cdot\nabla\right)\varphi-A^*\varphi+\nabla \pi =f &\text{in }Q, \\
			\nabla\cdot \varphi=0 &\text{in }Q, \\
			\varphi=0 & \text{on }\Sigma, \\
			\varphi(\cdot,T)=0 & \text{in }\Omega.
		\end{cases}
\end{equation}
Then we have the following result. 
\begin{lema} \label{lemma:23}
Assume $T_0>0$ and $T\in (0,T_0)$. Assume also $A_{i,j}\in \mathcal{M}_N(\mathbb{R})$, $B_{i,j}\in \mathbb{R}^N$,  $(i,j\in \{1,\ldots,m\}).$ 
Then for any $f\in L^2(Q)^{Nm}$, \eqref{Stokes} admits a unique solution 
$$
\varphi \in L^2(0,T;H^2(\Omega)^{Nm})\cap L^\infty(0,T;H^1(\Omega)^{Nm}) \cap H^1(0,T;\mathcal{H}), \quad \pi \in L^2(0,T; \left[H^1(\Omega)/\mathbb{R}\right]^m)
$$
and there exists a constant $C$ depending on $T_0$, on the geometry and on $A_{i,j}$, $B_{i,j}$ such that
\begin{equation}\label{regpara}
	\left\|\varphi \right\|_{L^2(0,T;H^2(\Omega)^{Nm})}+\left\|\varphi \right\|_{L^\infty(0,T;H^1(\Omega)^{Nm})}+\left\| \varphi \right\|_{H^1(0,T;L^2(\Omega)^{Nm})}
		\leq C\left\|f\right\|_{L^2(0,T;L^2(\Omega)^{Nm})}.
\end{equation}
Assume moreover that
$$
f\in L^2(0,T;\mathcal{V})
$$
then	
$$
\varphi \in L^2(0,T;H^3(\Omega)^{Nm})\cap L^\infty(0,T;H^2(\Omega)^{Nm})\cap H^1(0,T;\mathcal{V})
$$
and there exists a constant $C>0$ depending on $T_0$, on the geometry and on $A_{i,j}$, $B_{i,j}$ such that
\begin{equation}  \label{eq:h3v}
\left\|\varphi \right\|_{L^2\left(0,T;H^3\left(\Omega\right)^{Nm}\right)}
+\left\|\varphi \right\|_{L^\infty\left(0,T;H^2\left(\Omega\right)^{Nm}\right)}
+\left\|\varphi \right\|_{H^1\left(0,T;H^1(\Omega)^{Nm}\right)}\leq C\left\|f\right\|_{L^2\left(0,T;H^1(\Omega)^{Nm}\right)}.
\end{equation}
Assume moreover that 
\begin{equation}\label{hyp3}
f\in L^2\left(0,T;H^3(\Omega)^{Nm}\right)\cap H^1\left(0,T;\mathcal{V}\right), \quad f\left(\cdot,0\right)=0 \ \text{in} \ \Omega.
\end{equation}
Then, $\varphi\in L^2(0,T;H^5(\Omega)^{Nm})\cap H^1(0,T;H^3(\Omega)^{Nm})\cap H^2(0,T;\mathcal{V})$ 
and there exists a constant $C>0$ depending on $T_0$, on the geometry and on $A_{i,j}$, $B_{i,j}$ such that
$$			
\left\|\varphi \right\|_{L^2(0,T;H^5(\Omega)^{Nm})}
+\left\|\varphi\right\|_{H^1(0,T;H^3(\Omega)^{Nm})}
+\left\|\varphi\right\|_{H^2(0,T;H^1(\Omega)^{Nm})} 
\\
\leq C\left(\|f\|_{L^2(0,T;H^3(\Omega)^{Nm})}
+\|\partial_t f\|_{L^2(0,T;H^1(\Omega)^{Nm})}
\right).
$$
\end{lema}
\begin{proof}
We can write \eqref{Stokes} in the form
\begin{equation} \label{tt70}
    \begin{cases}
	    -\varphi'+\mathcal{A}^*\varphi=\mathcal{P}f \quad \text{in} \ (0,T),
	    \\
	    \varphi(T)=0.
    \end{cases}
\end{equation}

The first regularity result and the corresponding estimate \eqref{regpara} is thus a direct consequence from the fact that $-\mathcal{A}^*$ is the infinitesimal generator of an analytic semigroup.

Assume now that $f\in L^2(0,T;\mathcal{V})$.
We take the inner product of the first equation of \eqref{tt70} with $-\mathcal{A}^*\varphi'$ and we obtain
\begin{multline*}
\left\| \mathcal{A}_0^{1/2} \varphi' \right\|_{\mathcal{H}}^2-\frac{1}{2} \frac{d}{dt} \left\| \mathcal{A}^*\varphi \right\|_{\mathcal{H}}^2
= (\varphi',\mathcal{A}_1^* \varphi')_{\mathcal{H}}-\left( \mathcal{A}_0^{1/2} f, \mathcal{A}_0^{1/2}\varphi'\right)_{\mathcal{H}}
-\left( f, \mathcal{A}_1^*\varphi'\right)_{\mathcal{H}}
\\
\leq \frac{1}{2} \left\| \mathcal{A}_0^{1/2} \varphi' \right\|_{\mathcal{H}}^2
+C \left(\left\| \varphi' \right\|_{\mathcal{H}}^2 + \left\| f\right\|_{\mathcal{V}}^2\right).
\end{multline*}
Thus, using \eqref{regpara}, we deduce that
$$
\left\|\varphi \right\|_{H^1(0,T;\mathcal{V})} + 
\left\|\varphi \right\|_{L^\infty(0,T;\mathcal{D}(\mathcal{A}^*))}
\leq C \left\| f \right\|_{L^2(0,T;\mathcal{V})}.
$$
Going back to \eqref{Stokes}, this yields that
\begin{equation}\label{Stokes-k}
	\begin{cases}
		-\Delta \varphi+\nabla \pi =g:=f+\partial_t \varphi -\left(B^*\cdot\nabla\right)\varphi+A^*\varphi&\text{in }Q, \\
			\nabla\cdot \varphi=0 &\text{in }Q, \\
			\varphi=0 & \text{on }\Sigma, \\
		\end{cases}
\end{equation}
with $g\in L^2(0,T;H^1(\Omega)^{Nm})$. Using Proposition 2.2 p.33 in \cite{Temam}, we deduce that $\varphi\in L^2(0,T;H^3(\Omega)^{Nm})$ with the estimate 
\eqref{eq:h3v}.

Assume now \eqref{hyp3}. By taking the time derivative of \eqref{tt70}, we deduce that 
\begin{equation*} 
    \begin{cases}
	    -\left(\varphi'\right)'+\mathcal{A}^*\left(\varphi'\right)=f' \quad \text{in} \ (0,T),
	    \\
	    \left(\varphi'\right)(T)=0.
    \end{cases}
\end{equation*}
Thus, from the previous step, we deduce that
$$
\varphi \in H^1(0,T;H^3(\Omega)^{Nm})\cap W^{1,\infty}(0,T;H^2(\Omega)^{Nm})\cap H^2(0,T;\mathcal{V}).
$$
In particular, in \eqref{Stokes-k}, the right-hand side satisfies $g\in L^2(0,T;H^2(\Omega)^{Nm})$ and thus applying again 
Proposition 2.2 p.33 in \cite{Temam}, we deduce that $\varphi\in L^2(0,T;H^4(\Omega)^{Nm})$. This yields that $g\in L^2(0,T;H^3(\Omega)^{Nm})$, which allows us to conclude the proof of the lemma.
\end{proof}
	
\section{The Carleman estimate for the adjoint system}\label{S3}
	
As before, we denote by $C$ various positive constants which depends only on $\Omega$ and $\omega$ (they depend also in general on the choice of $\eta^0$ and $\omega_0$ but one can consider that $\eta^0$ as well as $\omega_0$ depend $\Omega$ and $\omega$). Without any lack of generality, we treat the case of dimension $N=2$. The same proof can be performed in the general case. 

Our aim is to estimate the following quantity associated with the solutions of the system \eqref{eq:systemadjoint2}:
	\begin{multline}\label{defI1}
		I(s,\varphi^{(i)}):=
		\iint_Q  e^{-2s\alpha}\left(\left(s\xi\right)^{-1} \left|\nabla^3 \Delta\varphi^{(i)}_1\right|^2
		+s\xi \left|\nabla^2 \Delta\varphi^{(i)}_1\right|^2
		+\left(s\xi\right)^3 \left|\nabla\Delta\varphi_1^{(i)}\right|^2
		+\left(s\xi\right)^5 \left|\Delta\varphi_1^{(i)}\right|^2
		\right)\diff{x}\diff{t}  
		\\
		+\iint_Qe^{-2s\alpha^*}(s\xi^*)^5 \left|\varphi^{(i)} \right|^2 \diff{x}\diff{t}.
	\end{multline}
Thus, our main result states as follows:
	\begin{teorema} \label{prop:0}
		There exists $C>0$ depending on the geometry such that for any $s\geq C$ and for any $\varphi_T\in \mathcal{H}$, the solution $\varphi=(\varphi^{(1)},...,\varphi^{(m)})$ of \eqref{eq:systemadjoint2} satisfies 
		\begin{equation} \label{eq:phiiphi11}
			\sum_{i=1}^m I(s,\varphi^{(i)})
			\leq C\iint_{\omega\times(0,T)}e^{-2s\alpha}s^{2^{m+3}-6}\xi^{2^{m+3}-6} \left| \varphi_1^{(1)} \right|^2  \diff{x} \diff{t}.
		\end{equation}
	\end{teorema}
	
In order to prove the above proposition, we first start by estimating each $I(s,\varphi^{(i)})$	$(i=1,\ldots,m)$ independently.
	\begin{proposicion} \label{prop:1}
		Let $\widehat{\omega}\subset\Omega$ be a nonempty open set such that $\omega_0\Subset\widehat\omega.$ Then, there exists a constant $C$ such that for any $s\geq C$
		and for any $\varphi_T\in \mathcal{H}$, the solution $\varphi$ of \eqref{eq:systemadjoint2} satisfies 
        \begin{multline} \label{eq:prop1}
				I(s,\varphi^{(i)})
				\leq C\left(
				\iint_{\widehat\omega\times(0,T)}e^{-2s\alpha} 
				s^5\xi^5 \left|\Delta\varphi_1^{(i)}\right|^2  \diff{x}\diff{t}
				\right. \\ \left.
				+\sum_{j=1}^{m} \iint_Q \left(e^{-2s\alpha^*}s^{9/2}(\xi^*)^{5} \left|\varphi^{(j)}\right|^2 
				+e^{-2s\alpha}s^{-2}\xi^{-2} \left(\left|\nabla^3 \Delta\varphi^{(j)}_1\right|^2+\left|\nabla^2 \Delta\varphi^{(j)}_1\right|^2\right)  \right) \diff{x}\diff{t}\right)
				\quad (1\leq i\leq m).
		\end{multline}
	\end{proposicion}
	\begin{proof}[Proof of Proposition \ref{prop:1}]
		First taking the divergence of \eqref{tt00}, we remark that
		$$
		\Delta \pi^{(i)}=0\quad\text{in }Q \quad (1\leq i\leq m).
		$$
		Thus, following the method introduced in \cite{Carreno-Gueye2014}, 
		we apply the operator $\nabla^2 \Delta$ to the first components of \eqref{tt00}, and we deduce 
		\begin{equation}\label{tt01}
				-\partial_t \nabla^2\Delta\varphi^{(i)}_1-\Delta \nabla^2\Delta\varphi^{(i)}_1 =-\sum_{j=1}^{i} B_{j,i}\cdot \nabla \left(\nabla^2\Delta\varphi^{(j)}_1\right)+\sum_{j=1}^{i+1}a_{j,i}\nabla^2\Delta\varphi^{(j)}_1, \quad (1\leq i \leq m).
		\end{equation}
		Applying Lemma \ref{lemma:22} to the above equations, we deduce that for $\lambda\geq\widehat{\lambda}_1$ and for $s\geq\widehat{s}_1$,
		\begin{multline} \label{eq:lemma5}
    	    \iint_Q  e^{-2s\alpha}\left(s^{-1}\xi^{-1} \left|\nabla^3 \Delta\varphi^{(i)}_1\right|^2
    	    +s\xi \left|\nabla^2 \Delta\varphi^{(i)}_1\right|^2\right)\diff{x}\diff{t}  
    	    \\
    	    \leq C\left(\left\|e^{-s\alpha^*}s^{-1/4}(\xi^*)^{-1/4+1/\ell} \nabla^2 \Delta\varphi^{(i)}_1\right\|^2_{L^2(\Sigma)^4} 
    	    +\left\|e^{-s\alpha^*}(s\xi^*)^{-1/4} \nabla^2 \Delta\varphi^{(i)}_1\right\|^{2}_{H^{\frac{1}{4},\frac{1}{2}}(\Sigma)^4}
    	    \right. \\ \left. 
    	    +\sum_{j=1}^{i} \iint_Q e^{-2s\alpha}s^{-2}\xi^{-2}\left|\nabla^3 \Delta\varphi^{(j)}_1\right|^2\diff{x}\diff{t}
    	    +\sum_{j=1}^{i+1} \iint_Q e^{-2s\alpha}s^{-2}\xi^{-2}\left|\nabla^2 \Delta\varphi^{(j)}_1\right|^2\diff{x}\diff{t}
    	    \right. \\ \left.
    	    +\iint_{\omega_0\times(0,T)}e^{-2s\alpha}s\xi\left|\nabla^2 \Delta\varphi^{(i)}_1\right|^2 \diff{x}\diff{t}\right).
		\end{multline}
		
		The rest of the proof is divided into several steps:
		\begin{itemize}
			\item In Step 1, we complete the left-hand side of \eqref{eq:lemma5} with weighted integrals of $\varphi^{(i)}$ in $Q$, and adding some local terms in the right-hand side.
			\item In Step 2, we obtain an upper bound of the boundary terms.
			\item Finally, in Step 3, we estimate the local terms that do not appear in \eqref{eq:prop1}.
		\end{itemize}
		
		\textbf{Step 1.} We apply Lemma \ref{lemma:21} with $u=\nabla\Delta\varphi_1^{(i)}$ and $r=1$: for any $s\geq C$ and $\lambda \geq C$,
		\begin{equation} \label{eq:lemma3r1}
			\iint_Q e^{-2s\alpha} s^3\xi^3 \left|\nabla\Delta\varphi_1^{(i)}\right|^2  \diff{x}\diff{t}
			\leq C\left( \iint_Q e^{-2s\alpha}s\xi \left|\nabla^2\Delta\varphi_1^{(i)}\right|^2 \diff{x}\diff{t}
			+\iint_{\omega_0\times(0,T)}e^{-2s\alpha}s^3\xi^3\left|\nabla\Delta\varphi_1^{(i)}\right|^2  \diff{x}\diff{t}
			\right).
		\end{equation}
		
		Then, we apply Lemma \ref{lemma:21} with $u=\Delta\varphi_1^{(i)}$ and $r=3$: for any $s\geq C$ and $\lambda \geq C$,
		\begin{equation} \label{eq:lemma3r3}
			\iint_Q e^{-2s\alpha} s^5\xi^5 \left|\Delta\varphi_1^{(i)}\right|^2  \diff{x}\diff{t}
			\leq 
			C\left(
			\iint_Q e^{-2s\alpha}s^3\xi^3\left|\nabla\Delta\varphi_1^{(i)}\right|^2 \diff{x}\diff{t}
			+\iint_{\omega_0\times(0,T)}e^{-2s\alpha} s^5\xi^5\left|\Delta\varphi_1^{(i)}\right|^2 \diff{x}\diff{t}
			\right).
		\end{equation}	
				
		Now, using the divergence condition of $\varphi^{(i)}$, we have
		$$
		\left|\partial_2 \varphi_2^{(i)}\right|=\left|\partial_1 \varphi_1^{(i)}\right| \leq \left|\nabla \varphi_1^{(i)}\right|.
		$$
		Then using the Poincar\'e inequality and the ellipticity of the Laplace operator with Dirichlet boundary conditions,
		we deduce the existence of a constant $C$ depending on $\Omega$ such that
		$$
		\int_{\Omega} \left|\varphi^{(i)}\right|^2  \diff{x}
		\leq C\int_{\Omega} \left|\Delta \varphi_1^{(i)}\right|^2  \diff{x}.
		$$
		Combining the above relation with \eqref{eq:lemma5}, \eqref{eq:lemma3r1} and \eqref{eq:lemma3r3}, we deduce
		that $I(s,\varphi^{(i)})$ defined by \eqref{defI1}
		satisfies 
		\begin{multline} \label{eq:sumI1}
			I(s,\varphi^{(i)})
			\leq C\left(
	        \iint_{\omega_0\times(0,T)}e^{-2s\alpha} 
        	\left(s\xi \left|\nabla^2 \Delta\varphi^{(i)}_1\right|^2
        	+s^3\xi^3 \left|\nabla\Delta\varphi_1^{(i)}\right|^2
        	+s^5\xi^5 \left|\Delta\varphi_1^{(i)}\right|^2 
	        \right) \diff{x}\diff{t}
	        \right. \\ \left.
	        +\left\|e^{-s\alpha^*}s^{-1/4}(\xi^*)^{-1/4+1/\ell} \nabla^2 \Delta\varphi^{(i)}_1\right\|^2_{L^2(\Sigma)^4}  
	        +\left\|e^{-s\alpha^*}(s\xi^*)^{-1/4} \nabla^2 \Delta\varphi^{(i)}_1\right\|^{2}_{H^{\frac{1}{4},\frac{1}{2}}(\Sigma)^4}
	        \right. \\ \left.
	        +\sum_{j=1}^{i} \iint_Qe^{-2s\alpha}s^{-2}\xi^{-2}\left|\nabla^3 \Delta\varphi^{(j)}_1\right|^2\diff{x}\diff{t}
	        +\sum_{j=1}^{i+1} \iint_Qe^{-2s\alpha}s^{-2}\xi^{-2}\left|\nabla^2 \Delta\varphi^{(j)}_1\right|^2\diff{x}\diff{t}
        	\right).
		\end{multline}
		
		\textbf{Step 2.} In this step, we get rid of the boundary terms in the right-hand side of \eqref{eq:sumI1}.
		To estimate the first term, we notice that since $\ell\geq 4$, $(\xi^*)^{-1/4+1/\ell}$ is bounded in $(0,T)$ (see \eqref{09:50}). Thus
		\begin{multline}\label{firstboundary}
			\left\|e^{-s\alpha^*}s^{-1/4}(\xi^*)^{-1/4+1/\ell} \nabla^2 \Delta\varphi^{(i)}_1\right\|^2_{L^2(\Sigma)^4}  
			\leq Cs^{-1/2}\left\| e^{-s\alpha^*}\nabla^2\Delta\varphi_1^{(i)} \right\|^2_{L^2(\Sigma)^4}  
			\\
			\leq Cs^{-1/2} \left(
			\left\|e^{-s\alpha^*}\nabla^2\Delta\varphi_1^{(i)} \right\|^2_{L^2(Q)^4}  
			+ \left\|e^{-s\alpha^*}(s\xi^*)^{1/2}\nabla^2\Delta\varphi_1^{(i)}\right\|_{L^2(Q)^4}
			\left\|e^{-s\alpha^*}(s\xi^*)^{-1/2}\nabla^3\Delta\varphi_1^{(i)}\right\|_{L^2(Q)^8}\right)  
			\\
			\leq C s^{-1/2}\iint_Qe^{-2s\alpha^*}\left(s\xi\left|\nabla^2\Delta\varphi_1^{(i)}\right|^2+s^{-1}\xi^{-1}\left|\nabla^3\Delta\varphi_1^{(i)}\right|^2\right)\diff{x}\diff{t}
			\leq Cs^{-1/2}I(s,\varphi^{(i)}).
		\end{multline}
		
		For the second boundary term, we will use a trace inequality (see, for instance, \cite[Theorem 2.1, p.9]{lionsnon2}) and by an interpolation argument 
		(see, for instance, \cite[Remark 9.5, pp.46-47]{lionsnon1}), we have
		\begin{multline} \label{eq:traceinequality}
			\left\| e^{-s\alpha^*} (\xi^*)^{-1/4}  \nabla^2\Delta\varphi_1^{(i)} \right\|^2_{H^{\frac{1}{4},\frac{1}{2}}(\Sigma)^4} 
			\leq C\left\| e^{-s\alpha^*} (\xi^*)^{-1/4}  \nabla^2\Delta\varphi_1^{(i)} \right\|^2_{L^2(0,T;H^1(\Omega)^4)\cap H^{1/2}(0,T;L^2(\Omega)^4)} 
			\\
			\leq C\left\| e^{-s\alpha^*} (\xi^*)^{-1/4}  \varphi_1^{(i)} \right\|^2_{L^2(0,T;H^5(\Omega)^4)\cap H^{1}(0,T;H^{3}(\Omega)^4)}. 
		\end{multline}
		Our goal is to estimate the last term.	First we consider the function 
		$$
		\theta_1(t):=s^{3/2}(\xi^*)^{3/2-1/\ell}e^{-s\alpha^*}.
		$$
		From \eqref{09:49} and \eqref{09:50}, for $s\geq C$,
		\begin{equation}\label{10:01}
			\left|\theta_{1}'\right| \leq C s^{5/2}(\xi^*)^{5/2} e^{-s\alpha^*}.
		\end{equation}
		Then, from \eqref{eq:adjointmatriz}, $(\theta_1\varphi,\theta_1\pi)$ is the solution of
		\begin{equation*} \label{eq:systemtheta}
		\begin{cases}
		    -\partial_t \left(\theta_1\varphi\right)-\Delta\left(\theta_1\varphi\right)+\left(B^*\cdot\nabla\right)\left(\theta_1\varphi\right)-A^*\left(\theta_1\varphi\right)+\nabla\left(\theta_1\pi\right)=-\theta_1'\varphi & \text{in }Q, \\
	    	\nabla\cdot\left(\theta_1\varphi\right)=0 & \text{in }Q, \\
	    	\left(\theta_1\varphi\right)=0 & \text{on }\Sigma, \\
	    	\left(\theta_1\varphi\right)(\cdot,T)=0 & \text{in }\Omega.
	    \end{cases}
		\end{equation*}{\tiny }
		Applying Lemma \ref{lemma:23} to the above system, we obtain
		\begin{equation} \label{eq:theta1leq}
			\left\|\theta_1\varphi\right\|^2_{L^2(0,T;H^2(\Omega)^{2m})\cap H^1(0,T;L^2(\Omega)^{2m})} \\
	        \leq C\left\|\theta_1'\varphi\right\|^2_{L^2(0,T;L^2(\Omega)^{2m})} 
	        \leq C\left\|e^{-s\alpha^*}s^{5/2}(\xi^*)^{5/2} \varphi\right\|^2_{L^2(0,T;L^2(\Omega)^{2m})}.
		\end{equation}

		Now, combining the above estimate with an interpolation inequality, we deduce 
		\begin{multline} \label{eq:estimateV}
	        \left\|e^{-s\alpha^*} s^2(\xi^*)^{2-1/({2\ell})} \varphi\right\|^2_{L^2(0,T;H^1(\Omega)^{2m})} 
        	\leq \left\|\theta_1\varphi\right\|_{L^2(0,T;H^2(\Omega)^{2m})}\left\|e^{-s\alpha^*}s^{5/2}(\xi^*)^{5/2}\varphi\right\|_{L^2(0,T;L^2(\Omega)^{2m})}  \\
         	\leq C\left\| e^{-s\alpha^*}s^{5/2}(\xi^*)^{5/2}\varphi\right\|^2_{L^2(0,T;L^2(\Omega)^{2m})}.
        \end{multline}		
		Second, we introduce
		$$
		\theta_2=s(\xi^*)^{1-3/{(2\ell)}}e^{-s\alpha^*}.
		$$
		From \eqref{09:49} and \eqref{09:50}, for $s\geq C$,
		\begin{equation}\label{15:57}
			\left|\theta_{2}'\right| \leq C s^{2}(\xi^*)^{2-1/(2\ell)} e^{-s\alpha^*}.
		\end{equation}
		Then, from \eqref{eq:adjointmatriz}, $(\theta_2\varphi,\theta_2\pi)$ is the solution of
		\begin{equation*} \label{eq:systemtheta2}
			\begin{cases}
				-\partial_t \left(\theta_2\varphi\right)-\Delta\left(\theta_2\varphi\right)+\left(B^*\cdot\nabla\right)\left(\theta_2\varphi\right)-A^*\left(\theta_2\varphi\right)+\nabla\left(\theta_2\pi\right)=-\theta_2'\varphi & \text{in }Q, \\
				\nabla\cdot\left(\theta_2\varphi\right)=0 & \text{in }Q, \\
				\left(\theta_2\varphi\right)=0 & \text{on }\Sigma, \\
				\left(\theta_2\varphi\right)(\cdot,T)=0 & \text{in }\Omega.
			\end{cases}
		\end{equation*}
		Since $-\theta_2'\varphi\in L^2(0,T;\mathcal{V})$, we can apply 
		Lemma \ref{lemma:23} to the above system and we deduce 
		\begin{equation*}
    	    \left\|\theta_2\varphi\right\|^2_{L^2(0,T;H^3(\Omega)^{2m})\cap H^1(0,T;H^1(\Omega)^{2m})}
	        \leq C\left\|\theta_2'\varphi\right\|^2_{L^2(0,T;H^1(\Omega)^{2m})}.
        \end{equation*}
        Using \eqref{15:57} and combining the above result with \eqref{eq:estimateV}, we obtain
			\begin{equation}\label{16:16}
			\left\| \theta_2\varphi \right\|^2_{L^2(0,T;H^3(\Omega)^2)\cap H^1(0,T;H^1(\Omega)^{2m})}
			\leq C\left\|s^{5/2}(\xi^*)^{5/2} e^{-s\alpha^*}\varphi\right\|^2_{L^2(0,T;L^2(\Omega)^{2m})}.
		\end{equation}
		Finally, we introduce
		$$
		\theta_3=(\xi^*)^{-5/{(2\ell)}}e^{-s\alpha^*}.
		$$
		From \eqref{09:49} and \eqref{09:50}, for $s\geq C$,
		\begin{equation}\label{18:11}
			\left|\theta_{3}'\right| \leq C \theta_2, \quad
			\left|\theta_{3}''\right| \leq C s^2(\xi^*)^{2-1/({2\ell})} e^{-s\alpha^*}.
		\end{equation}
		Then, from \eqref{eq:adjointmatriz}, $(\theta_3\varphi,\theta_3\pi)$ is the solution of
		\begin{equation*} \label{eq:systemtheta3}
			\begin{cases}
				-\partial_t \left(\theta_3\varphi\right)-\Delta\left(\theta_3\varphi\right)+\left(B^*\cdot\nabla\right)\left(\theta_3\varphi\right)-A^*\left(\theta_3\varphi\right)+\nabla\left(\theta_3\pi\right)=-\theta_3'\varphi & \text{in }Q, \\ 
				\nabla\cdot\left(\theta_3\varphi\right)=0 & \text{in }Q, \\
				\left(\theta_3\varphi\right)=0 & \text{on }\Sigma, \\
				\left(\theta_3\varphi\right)(\cdot,T)=0 & \text{in }\Omega.
			\end{cases}
		\end{equation*}
		Since $-\theta_3'\varphi$ satisfies \eqref{hyp3}, we can apply \cref{lemma:23} to the above system and we deduce	
		\begin{multline}\label{eq:estimateH5-1}
			\left\|\theta_3\varphi\right\|^2_{L^2(0,T;H^5(\Omega)^{2m})\cap H^1(0,T;H^3(\Omega)^{2m})\cap H^2(0,T;H^1(\Omega)^{2m})} \\
        	\leq C\left(\left\|\theta_3'\varphi\right\|^2_{L^2(0,T;H^3(\Omega)^{2m})}+\left\|\partial_t(\theta_3'\varphi)\right\|^2_{L^2(0,T;H^1(\Omega)^{2m)}}
        	\right).
		\end{multline}
		
		Using \eqref{18:11} and applying \eqref{16:16} and \eqref{eq:estimateV}, we deduce from \eqref{eq:estimateH5-1}
		\begin{equation} \label{eq:estimateH5}
			\left\|\theta_3\varphi\right\|^2_{L^2(0,T;H^5(\Omega)^{2m})\cap H^1(0,T;H^3(\Omega)^{2m})\cap H^2(0,T;H^1(\Omega)^{2m})} 
			\leq C\left\|e^{-s\alpha^*} s^{5/2}(\xi^*)^{5/2} \varphi\right\|^2_{L^2(0,T;L^2(\Omega)^{2m})}.
		\end{equation}
		Using the above estimate, \eqref{09:49} and that $\ell\geq 14$, we have
		$$
		\left\|e^{-s\alpha^*} (\xi^*)^{-1/4} \varphi\right\|^2_{L^2(0,T;H^5(\Omega)^{2m})\cap H^1(0,T;H^{3}(\Omega)^{2m})} 
		\leq C\left\|e^{-s\alpha^*}s^{5/2}(\xi^*)^{5/2} \varphi\right\|^2_{L^2(0,T;L^2(\Omega)^{2m})}.
		$$
		Combining the above estimate and the trace inequality \eqref{eq:traceinequality}, we deduce that
			$$
			\left\|e^{-s\alpha^*}(s\xi^*)^{-1/4} \nabla^2 \Delta\varphi^{(i)}_1\right\|^{2}_{H^{\frac{1}{4},\frac{1}{2}}(\Sigma)^4} 
			\leq Cs^{-1/2}\left\|e^{-s\alpha^*}s^{5/2}(\xi^*)^{5/2} \varphi\right\|^2_{L^2(0,T;L^2(\Omega)^{2m})}.
			$$
		Putting together \eqref{eq:sumI1}, the above relation and \eqref{firstboundary}, we deduce at this step the following inequality: 
		\begin{multline} \label{eq:sumI2}
			I(s,\varphi^{(i)})
			\leq C\left(
            \iint_{\omega_0\times(0,T)}e^{-2s\alpha} \left(s\xi \left|\nabla^2 \Delta\varphi^{(i)}_1\right|^2
            	+s^3\xi^3\left|\nabla\Delta\varphi_1^{(i)}\right|^2
            	+s^5\xi^5\left|\Delta\varphi_1^{(i)}\right|^2 \right) \diff{x}\diff{t} 
            	\right. \\ \left.
            	+\sum_{j=1}^{m} \iint_Q \left(e^{-2s\alpha^*}s^{9/2}(\xi^*)^{5} \left|\varphi^{(j)}\right|^2 
            	+e^{-2s\alpha}(s\xi)^{-2}\left( \left|\nabla^3\Delta\varphi^{(j)}_1\right|^2
            	+\left|\nabla^2\Delta\varphi^{(j)}_1\right|^2 \right)
            	\right)  \diff{x}\diff{t} \right).
		\end{multline}
		\textbf{Step 3.} To estimate the local terms, we proceed in a standard way: we consider $\omega_1$ an open subset satisfying 
		$\omega_0\Subset\omega_1\Subset\widehat{\omega}$ and 
		$$
		\eta_1\in C_c^{2}(\omega_1), \quad \eta_1\equiv 1 \ \text{in} \ \omega_0, \quad \eta_1\geq 0.
		$$ 
		Then, an integration by parts gives 
		\begin{multline*}
			\iint_{\omega_0\times(0,T)} e^{-2s\alpha} s\xi \left(\frac{\partial^2}{\partial x_k\partial x_q} \Delta\varphi_1^{(i)}\right)^2  \diff{x}\diff{t}
			\leq 
			\iint_{\omega_1\times(0,T)} \eta_1 e^{-2s\alpha} s\xi \left(\frac{\partial^2}{\partial x_k\partial x_q} \Delta\varphi_1^{(i)}\right)^2  \diff{x}\diff{t}
			\\
			=-\iint_{\omega_1\times(0,T)} \frac{\partial}{\partial x_k} \left(\eta_1 e^{-2s\alpha} s\xi\right) \frac{\partial^2}{\partial x_k\partial x_q} \Delta\varphi_1^{(i)} 
			\frac{\partial}{\partial x_q} \Delta\varphi_1^{(i)}  \diff{x}\diff{t}
			\\
			-\iint_{\omega_1\times(0,T)} \eta_1 e^{-2s\alpha} s\xi \frac{\partial^3}{\partial x_k^2\partial x_q} \Delta\varphi_1^{(i)} 
			\frac{\partial}{\partial x_q} \Delta\varphi_1^{(i)}  \diff{x}\diff{t}.
		\end{multline*}
		Using \eqref{11:01} and Young's inequality, we deduce from the above relation that there exists $C>0$ such that for all $\varepsilon>0$,
		\begin{multline}\label{11:25}
			\iint_{\omega_0\times(0,T)}e^{-2s\alpha} 	s\xi \left|\nabla^2 \Delta\varphi^{(i)}_1\right|^2  \diff{x}\diff{t}
			\leq \varepsilon \iint_Q  e^{-2s\alpha}\left(s^{-1}\xi^{-1} \left|\nabla^3 \Delta\varphi^{(i)}_1\right|^2
			+s\xi \left|\nabla^2 \Delta\varphi^{(i)}_1\right|^2\right)\diff{x}\diff{t}  
			\\
			+\frac{C}{\varepsilon}\iint_{\omega_1\times(0,T)}e^{-2s\alpha} s^3\xi^3\left|\nabla\Delta\varphi_1^{(i)}\right|^2 \diff{x}\diff{t}.
		\end{multline}
		
		Now we estimate, in an analogous way, the local term associated with $\nabla\Delta\varphi_1^{(i)}$: we consider
		$$
		\eta_2\in C_c^{2}(\widehat\omega), \quad \eta_2\equiv 1 \ \text{in} \ \omega_1, \quad \eta_2\geq 0.
		$$ 
		Then, integrating by parts, we obtain 
		\begin{multline*}
			\iint_{\omega_1\times(0,T)} e^{-2s\alpha} s^3\xi^3 \left(\frac{\partial}{\partial x_k} \Delta\varphi_1^{(i)}\right)^2  \diff{x}\diff{t}
			\leq 
			\iint_{\widehat\omega\times(0,T)} \eta_2 e^{-2s\alpha} s^3\xi^3 \left(\frac{\partial}{\partial x_k} \Delta\varphi_1^{(i)}\right)^2  \diff{x}\diff{t}
			\\
			=-\iint_{\widehat\omega\times(0,T)} \frac{\partial}{\partial x_k} \left(\eta_2 e^{-2s\alpha} s^3\xi^3\right) \frac{\partial}{\partial x_k} \Delta\varphi_1^{(i)} 
			\Delta\varphi_1^{(i)}  \diff{x}\diff{t}
			\\
			-\iint_{\widehat\omega\times(0,T)} \eta_2 e^{-2s\alpha} s^3\xi^3 \frac{\partial^2}{\partial x_k^2} \Delta\varphi_1^{(i)} \Delta\varphi_1^{(i)}  \diff{x}\diff{t}.
		\end{multline*}
		
		Using \eqref{11:01} and the Young's inequality, we deduce from the above relation that there exists $C>0$ such that for all $\varepsilon>0$,
		\begin{multline*}
			\iint_{\omega_1\times(0,T)}e^{-2s\alpha} 	s^3\xi^3 \left|\nabla \Delta\varphi^{(i)}_1\right|^2  \diff{x}\diff{t}
			\leq \varepsilon \iint_Q  e^{-2s\alpha}\left(s^{3}\xi^{3} \left|\nabla \Delta\varphi^{(i)}_1\right|^2
			+s\xi \left|\nabla^2 \Delta\varphi^{(i)}_1\right|^2
			\right)\diff{x}\diff{t}  
			\\
			+\frac{C}{\varepsilon}\iint_{\widehat{\omega}\times(0,T)}e^{-2s\alpha} s^5\xi^5 \left|\Delta\varphi_1^{(i)}\right|^2 \diff{x}\diff{t}.
		\end{multline*}
		The above estimate, together with \eqref{eq:sumI2} and \eqref{11:25} implies \eqref{eq:prop1}. This concludes the proof of Proposition \ref{prop:1}.
	\end{proof}
	
	We are now in a position to prove \cref{prop:0}.
	\begin{proof}[Proof of \cref{prop:0}]
		Summing \eqref{eq:prop1} for $i=1,\ldots,m$, and taking $s\geq C$ for a constant $C$ large enough, we deduce that
		\begin{equation} \label{sumavarphii2}
			\sum_{i=1}^m I(s,\varphi^{(i)})
			\leq C
			\sum_{i=1}^m\iint_{\widehat\omega\times(0,T)}e^{-2s\alpha} 
			s^5\xi^5 \left|\Delta\varphi_1^{(i)}\right|^2  \diff{x}\diff{t}.
		\end{equation}
		In order to get rid of the local terms in the right-hand side (except the term corresponding to $i=1$), we introduce a sequence of open sets
		$\mathcal{O}_i$, $(0\leq i\leq m)$ such that
		$$
		\widehat\omega=:\mathcal{O}_0\Subset\mathcal{O}_1\Subset...\Subset\mathcal{O}_i\Subset...\Subset\mathcal{O}_m\Subset\omega
		$$
		and functions
		$$
		\zeta_i\in C_c^2(\mathcal{O}_i) \quad \text{such that}\quad \zeta_i\equiv 1 \ \text{in}\ \mathcal{O}_{i-1}, \quad \zeta_i\geq 0 \quad (1\leq i\leq m).
		$$
		Then, we consider the equation $m-1$ of \eqref{tt00}, we apply the Laplace operator on the first component of this equation and we multiply it by $\zeta_1$:
		\begin{equation} \label{eq:mm-1}
			-\zeta_1\partial_t \Delta\varphi^{(m-1)}_1-\zeta_1\Delta^2 \varphi^{(m-1)}_1 +\zeta_1\sum_{j=1}^{m-1}B_{j,m-1}\cdot \nabla \Delta\varphi^{(j)}_1-\zeta_1\sum_{j=1}^{m-1}a_{j,m-1}\Delta\varphi^{(j)}_1=\zeta_1a_{m,m-1}\Delta\varphi^{(m)}_1
		\end{equation}
		
		Then, using the above equation, we deduce
\begin{equation}\label{16:51}
			\iint_{\widehat\omega\times(0,T)}e^{-2s\alpha} 
	        s^5\xi^5 \left|\Delta\varphi_1^{(m)}\right|^2 \diff{x}\diff{t}
        	\leq 
        	\iint_{\mathcal{O}_1\times(0,T)} \zeta_1 e^{-2s\alpha} s^5\xi^5 \left|\Delta\varphi_1^{(m)}\right|^2 \diff{x}\diff{t}
	=\sum_{k=1}^4 J_k
\end{equation}
with
\begin{equation}\label{16:51-J1}
J_1:=-\frac{1}{a_{m,m-1}}\iint_{\mathcal{O}_1\times(0,T)} \zeta_1 
	        e^{-2s\alpha} s^5\xi^5 \left(\Delta\varphi_1^{(m)}\right)\left(\partial_t\Delta\varphi^{(m-1)}_1 \right)\diff{x}\diff{t},
\end{equation}
\begin{equation}\label{16:51-J2}
J_2:=-\frac{1}{a_{m,m-1}}\iint_{\mathcal{O}_1\times(0,T)} \zeta_1 
	        e^{-2s\alpha} s^5\xi^5 \left(\Delta\varphi_1^{(m)}\right)\left(\Delta^2\varphi_1^{(m-1)} \right)\diff{x}\diff{t},
\end{equation}
\begin{equation}\label{16:51-J3}
J_3:=\frac{1}{a_{m,m-1}}\iint_{\mathcal{O}_1\times(0,T)} \zeta_1 
	        e^{-2s\alpha} s^5\xi^5 \left(\Delta\varphi_1^{(m)}\right)\left(\sum_{j=1}^{m-1}B_{j,m-1}\cdot \nabla \Delta\varphi^{(j)}_1 \right)\diff{x}\diff{t},
\end{equation}
\begin{equation}\label{16:51-J4}
J_4:=-\frac{1}{a_{m,m-1}}\iint_{\mathcal{O}_1\times(0,T)} \zeta_1 
	        e^{-2s\alpha} s^5\xi^5 \left(\Delta\varphi_1^{(m)}\right)\left(\sum_{j=1}^{m-1}a_{j,m-1}\Delta\varphi^{(j)}_1 \right)\diff{x}\diff{t}.
\end{equation}
	
Let us start by estimating the term $J_1$. Integrating by parts and using
 \begin{equation*}
\zeta_1\partial_t \Delta\varphi^{(m)}_1=-\zeta_1\Delta^2\varphi^{(m)}_1 
+\zeta_1\sum_{j=1}^{m} B_{j,m}\cdot \nabla\Delta\varphi^{(j)}_1-\zeta_1\sum_{j=1}^{m}a_{j,m}\Delta\varphi^{(j)}_1,
\end{equation*}
we obtain
\begin{equation}\label{23:12}
J_1=\sum_{k=1}^4 J_{1,k}
\end{equation}
with
\begin{equation}\label{23:12-J11}
J_{1,1}:=\frac{1}{a_{m,m-1}} \iint_{\mathcal{O}_1\times(0,T)} \zeta_1 
    	    \partial_t\left(e^{-2s\alpha} s^5\xi^5\right) \Delta\varphi_1^{(m)} \Delta\varphi^{(m-1)}_1 \diff{x}\diff{t},
\end{equation}
\begin{equation}\label{23:12-J12}
J_{1,2}:=\frac{-1}{a_{m,m-1}} \iint_{\mathcal{O}_1\times(0,T)} \zeta_1 
        	e^{-2s\alpha} s^5\xi^5 \Delta^2\varphi^{(m)}_1 \Delta\varphi^{(m-1)}_1 \diff{x}\diff{t},
\end{equation}
\begin{equation}\label{23:12-J13}
J_{1,3}:=\frac{1}{a_{m,m-1}}  \iint_{\mathcal{O}_1\times(0,T)} \zeta_1 
        	e^{-2s\alpha} s^5\xi^5  \left(\sum_{j=1}^{m} B_{j,m}\cdot \nabla\Delta\varphi^{(j)}_1\right) \Delta\varphi^{(m-1)}_1 \diff{x}\diff{t}, 
\end{equation}
\begin{equation}\label{23:12-J14}
J_{1,4}:=
\frac{1}{a_{m,m-1}}  \iint_{\mathcal{O}_1\times(0,T)} \zeta_1 
        	e^{-2s\alpha} s^5\xi^5  \left(\sum_{j=1}^{m}a_{j,m}\Delta\varphi^{(j)}_1\right) \Delta\varphi^{(m-1)}_1 \diff{x}\diff{t}.
\end{equation}
Using the estimate \eqref{09:48} and Young's inequality, we deduce the existence of $C$ such that for any $s\geq C$ and for any $\varepsilon >0,$
\begin{equation} \label{eq:J1}
        	\left| J_{1,1} \right| \leq \varepsilon \iint_Qe^{-2s\alpha}s^5\xi^5 \left|\Delta\varphi_1^{(m)}\right|^2 \diff{x}\diff{t}
        	+ \frac{C}{\varepsilon}\iint_{\mathcal{O}_1\times(0,T)} e^{-2s\alpha}s^7\xi^{7+2/\ell} \left|\Delta\varphi_1^{(m-1)}\right|^2 \diff{x}\diff{t},
\end{equation}
\begin{equation} \label{eq:J2}
        \left| J_{1,2} \right| \leq \varepsilon \iint_Q e^{-2s\alpha}s\xi \left| \Delta^2\varphi_1^{(m)}\right|^2 \diff{x}\diff{t}
        	+ \frac{C}{\varepsilon}\iint_{\mathcal{O}_1\times(0,T)} e^{-2s\alpha}s^9\xi^{9} \left|\Delta\varphi_1^{(m-1)}\right|^2 \diff{x}\diff{t},
\end{equation}
\begin{equation} \label{eq:J2.5}
        	\left| J_{1,3} \right| \leq        	\varepsilon \sum_{j=1}^{m}\iint_Qe^{-2s\alpha}s^3\xi^3 \left| \nabla\Delta\varphi_1^{(j)} \right|^2\diff{x}\diff{t}
        	+\frac{C}{\varepsilon}\iint_{\mathcal{O}_1\times(0,T)}e^{-2s\alpha}s^7\xi^7\left|\Delta\varphi_1^{(m-1)}\right|^2\diff{x}\diff{t},
\end{equation}
and
\begin{equation} \label{eq:J3}		
        \left| J_{1,4} \right| \leq
        	\varepsilon \sum_{j=1}^{m}\iint_Qe^{-2s\alpha}s^5\xi^5 \left| \Delta\varphi_1^{(j)} \right|^2\diff{x}\diff{t}
        	+\frac{C}{\varepsilon}\iint_{\mathcal{O}_1\times(0,T)}e^{-2s\alpha}s^5\xi^5\left|\Delta\varphi_1^{(m-1)}\right|^2\diff{x}\diff{t}.
        \end{equation}
Combining \eqref{eq:J1}, \eqref{eq:J2}, \eqref{eq:J2.5} and \eqref{eq:J3} we obtain
        \begin{align} \label{eq:estimateJ1}
        	\left| J_1\right| \leq \varepsilon\sum_{j=1}^m I(s,\varphi^{(j)})+\frac{C}{\varepsilon}\iint_{\mathcal{O}_1\times(0,T)} e^{-2s\alpha}s^9\xi^{9} \left|\Delta\varphi_1^{(m-1)}\right|^2 \diff{x}\diff{t}.
        \end{align}

        To estimate $J_2$, note that by integrating by parts, we find
        \begin{multline}\label{23:13}
        	\iint_{\mathcal{O}_1\times(0,T)} \zeta_1 
        	e^{-2s\alpha} s^5\xi^5 \left(\Delta\varphi_1^{(m)}\right)\left(\Delta^2\varphi_1^{(m-1)}\right)  \diff{x}\diff{t}
        	=\iint_{\mathcal{O}_1\times(0,T)}  
	        \Delta\left(\zeta_1 e^{-2s\alpha} s^5\xi^5\right) \Delta\varphi_1^{(m)} \Delta\varphi_1^{(m-1)}  \diff{x}\diff{t}
        	\\
        	+\iint_{\mathcal{O}_1\times(0,T)} \zeta_1 
        	e^{-2s\alpha} s^5\xi^5 \Delta^2\varphi_1^{(m)} \Delta\varphi_1^{(m-1)}  \diff{x}\diff{t}
	        \\
	        +\iint_{\mathcal{O}_1\times(0,T)} 
        	2\nabla\left(\zeta_1 e^{-2s\alpha} s^5\xi^5\right)\cdot \nabla \Delta\varphi_1^{(m)} \Delta\varphi_1^{(m-1)}  \diff{x}\diff{t}.
        \end{multline}
        Using \eqref{11:01} and Young's inequality, we deduce the existence of $C$ such that for any $s\geq C$ and for any $\varepsilon >0,$
        \begin{multline} \label{eq:J4}
        	\left| \iint_{\mathcal{O}_1\times(0,T)}  
        	\Delta\left(\zeta_1 e^{-2s\alpha} s^5\xi^5\right) \Delta\varphi_1^{(m)} \Delta\varphi_1^{(m-1)}  \diff{x}\diff{t}
        	\right|
        	\\
        	\leq \varepsilon \iint_Qe^{-2s\alpha}s^5\xi^5 \left|\Delta\varphi_1^{(m)}\right|^2 \diff{x}\diff{t}
        	+ \frac{C}{\varepsilon}\iint_{\mathcal{O}_1\times(0,T)} e^{-2s\alpha}s^9\xi^{9} \left|\Delta\varphi_1^{(m-1)}\right|^2 \diff{x}\diff{t}.
        \end{multline}		
        Again, using \eqref{11:01} and Young's inequality, we deduce the existence of $C$ such that for any $s\geq C$ and for any $\varepsilon >0,$
        \begin{multline} \label{eq:J5}
        	\left| 
        	\iint_{\mathcal{O}_1\times(0,T)} 
        	2\nabla\left(\zeta_1 e^{-2s\alpha} s^5\xi^5\right)\cdot \nabla \Delta\varphi_1^{(m)} \Delta\varphi_1^{(m-1)}  \diff{x}\diff{t}
        	\right|
        	\\
        	\leq \varepsilon \iint_Q e^{-2s\alpha}s^3\xi^3 \left|\nabla \Delta\varphi_1^{(m)}\right|^2 \diff{x}\diff{t}
        	+ \frac{C}{\varepsilon}\iint_{\mathcal{O}_1\times(0,T)} e^{-2s\alpha}s^9\xi^{9} \left|\Delta\varphi_1^{(m-1)}\right|^2 \diff{x}\diff{t}.
        \end{multline}		
        Combining \eqref{eq:J2}, \eqref{eq:J4} and \eqref{eq:J5} we get 
        \begin{align} \label{eq:estimateJ2}
        	\left| J_2 \right| \leq \varepsilon I(s,\varphi^{(m)})+\frac{C}{\varepsilon}\iint_{\mathcal{O}_1\times(0,T)} e^{-2s\alpha}s^9\xi^{9} \left|\Delta\varphi_1^{(m-1)}\right|^2 \diff{x}\diff{t}.
        \end{align}
We proceed similarly for $J_3$ (see \eqref{16:51-J3}), and after an integration by parts, we deduce the existence of $C$ such that for any $s\geq C$ and for any $\varepsilon>0$,
\begin{multline} \label{eq:J7}		
\left| J_3\right|=
\left| 
\frac{-1}{a_{m,m-1}}\sum_{j=1}^{m-1}B_{j,m-1}\cdot
\iint_{\mathcal{O}_1\times(0,T)} 
\left[
\nabla \left(\zeta_1 e^{-2s\alpha} s^5\xi^5\right) \Delta\varphi_1^{(m)}
+ \left(\zeta_1 e^{-2s\alpha} s^5\xi^5\right) \nabla\Delta\varphi_1^{(m)}
\right]
 \Delta\varphi^{(j)}_1 \diff{x}\diff{t}
\right|
        	\\
        	\leq 
        	\varepsilon  \iint_Qe^{-2s\alpha} \left( s^3 \xi^3 \left|\nabla\Delta\varphi_1^{(m)} \right|^2 
        	+ s^5 \xi^5 \left|\Delta\varphi_1^{(m)} \right|^2 \right)\diff{x}\diff{t} 
        	+ \frac{C}{\varepsilon} \sum_{j=1}^{m-1}\iint_{\mathcal{O}_1\times(0,T)} e^{-2s\alpha}s^7\xi^7\left|\Delta\varphi_1^{(j)}\right|^2 \diff{x}\diff{t}.
\end{multline}
        
        Finally, for $J_4$ (see \eqref{16:51-J3}), using Young's inequality, we deduce the existence of $C$ such that for any $s\geq C$ and for any $\varepsilon>0$
        \begin{equation} \label{eq:J6}		
        	\left| J_4 \right| 
           	\leq 
          	\varepsilon \iint_Qe^{-2s\alpha}s^5 \xi^5 \left|\Delta\varphi_1^{(m)} \right|^2 \diff{x}\diff{t}
        	+ \frac{C}{\varepsilon} \sum_{j=1}^{m-1}\iint_{\mathcal{O}_1\times(0,T)} e^{-2s\alpha}s^5\xi^5\left|\Delta\varphi_1^{(j)}\right|^2 \diff{x}\diff{t},
        \end{equation}
        for every $s\geq C$ and any $\varepsilon>0.$

        The combination of \eqref{16:51} with \eqref{eq:estimateJ1}, \eqref{eq:estimateJ2}, \eqref{eq:J7} and \eqref{eq:J6} yields the existence of a constant $C$ such that for any $s\geq C$ and for any $\varepsilon>0$,
        \begin{multline*}
        	\iint_{\widehat\omega\times(0,T)}e^{-2s\alpha} 
    	    s^5\xi^5 \left|\Delta\varphi_1^{(m)}\right|^2 \diff{x}\diff{t}
        	\leq \varepsilon \sum_{j=1}^{m} I(s,\varphi^{(j)})
        	\\
        	+ \frac{C}{\varepsilon} \iint_{\mathcal{O}_1\times(0,T)} e^{-2s\alpha} \left( s^9\xi^{9} \left|\Delta\varphi_1^{(m-1)}\right|^2 \diff{x}\diff{t}
        	+ \sum_{j=1}^{m-2}s^7\xi^{7} \left|\Delta\varphi_1^{(j)}\right|^2 \right) \diff{x}\diff{t} 
        	\\
        	\leq \varepsilon \sum_{j=1}^{m} I(s,\varphi^{(j)})
        	+ \frac{C}{\varepsilon}\sum_{j=1}^{m-1}\iint_{\mathcal{O}_1\times(0,T)} e^{-2s\alpha}s^9\xi^{9} \left|\Delta\varphi_1^{(j)}\right|^2 \diff{x}\diff{t}.
        \end{multline*}
        Analogously it can be proved that for $\Delta\varphi_1^{(m-1)}$, there exists a constant $C$ such that for any $s\geq C$ and for any $\varepsilon>0$,
        \begin{equation*}
        	\iint_{\mathcal{O}_1\times(0,T)}e^{-2s\alpha} 
        	s^{9}\xi^{9} \left|\Delta\varphi_1^{(m-1)}\right|^2 \diff{x}\diff{t}
        	\leq \varepsilon \sum_{j=1}^{m} I(s,\varphi^{(j)})
	        + \frac{C}{\varepsilon}\sum_{j=1}^{m-2}\iint_{\mathcal{O}_2\times(0,T)} e^{-2s\alpha}s^{17}\xi^{17} \left|\Delta\varphi_1^{(j)}\right|^2 \diff{x}\diff{t}.
        \end{equation*}
        
        Iterating the argument, we can estimate all the local terms and we deduce from \eqref{sumavarphii2} that
        \begin{equation} \label{eq:finiteiteration}
        	\sum_{i=1}^m I(s,\varphi^{(i)})
            \leq C\iint_{\mathcal{O}_m\times(0,T)}e^{-2s\alpha} 
        	s^{2^{(m+1)}+1}\xi^{2^{(m+1)}+1} \left|\Delta\varphi_1^{(1)}\right|^2 \diff{x}\diff{t}.
        \end{equation}

		Finally, we estimate the above local term in terms of $\varphi_1^{(1)}$. In order to do this, we consider
		$\widetilde\omega$ an open subset satisfying $\mathcal{O}_m\Subset\widetilde\omega\Subset \omega$ and 
		$$
		\widetilde\zeta\in C_c^2(\widetilde\omega) \quad \text{such that}\quad \widetilde\zeta\equiv 1 \ \text{in}\ \mathcal{O}_{m}, \quad \widetilde\zeta\geq 0.
		$$
		Then by integrating by parts, we obtain
		\begin{multline*}
			\iint_{\mathcal{O}_m\times(0,T)}e^{-2s\alpha} 
			s^{2^{(m+1)}+1}\xi^{2^{(m+1)}+1} \left|\Delta\varphi_1^{(1)}\right|^2  \diff{x}\diff{t}
			\leq 
			\iint_{\widetilde\omega\times(0,T)}\widetilde\zeta e^{-2s\alpha} 
			s^{2^{(m+1)}+1}\xi^{2^{(m+1)}+1} \left|\Delta\varphi_1^{(1)}\right|^2  \diff{x}\diff{t}
			\\
			=-\iint_{\widetilde\omega\times(0,T)}
			\nabla\left(\widetilde\zeta e^{-2s\alpha} s^{2^{(m+1)}+1}\xi^{2^{(m+1)}+1}\right) \Delta\varphi_1^{(1)}\cdot \nabla\varphi_1^{(1)}   \diff{x}\diff{t}
			\\
			-\iint_{\widetilde\omega\times(0,T)}
			\widetilde\zeta e^{-2s\alpha} s^{2^{(m+1)}+1}\xi^{2^{(m+1)}+1} \Delta\nabla\varphi_1^{(1)}\cdot \nabla\varphi_1^{(1)}   \diff{x}\diff{t}.
		\end{multline*}
		Considering \eqref{11:01}, and using Young's inequality, we deduce the existence of $C$ such that for any $s\geq C$ and for any $\varepsilon >0,$
		\begin{multline}  \label{eq:deltavarphi1nablavarphi1}
			\iint_{\mathcal{O}_m\times(0,T)}e^{-2s\alpha} 
			s^{2^{(m+1)}+1}\xi^{2^{(m+1)}+1} \left|\Delta\varphi_1^{(1)}\right|^2  \diff{x}\diff{t}
			\leq \varepsilon \iint_Q e^{-2s\alpha} \left(s^5\xi^5 \left| \Delta\varphi_1^{(1)}\right|^2+s^3\xi^3 \left|\nabla \Delta\varphi_1^{(1)}\right|^2\right) \diff{x}\diff{t}
			\\
			+ \frac{C}{\varepsilon}\iint_{\widetilde\omega\times(0,T)} e^{-2s\alpha}s^{2^{(m+2)}-1}\xi^{2^{(m+2)}-1} \left|\nabla \varphi_1^{(1)}\right|^2 \diff{x}\diff{t}.
		\end{multline}		
		Then, we consider
		$$
		\zeta\in C_c^2(\omega) \quad \text{such that}\quad \zeta\equiv 1 \ \text{in}\ \widetilde\omega, \quad \zeta\geq 0
		$$
		and we integrate by parts:
		\begin{multline*}
			\iint_{\widetilde\omega\times(0,T)} e^{-2s\alpha}s^{2^{(m+2)}-1}\xi^{2^{(m+2)}-1} \left|\nabla \varphi_1^{(1)}\right|^2 \diff{x}\diff{t}
			\leq
			\iint_{\omega\times(0,T)} \zeta e^{-2s\alpha}s^{2^{(m+2)}-1}\xi^{2^{(m+2)}-1} \left|\nabla \varphi_1^{(1)}\right|^2 \diff{x}\diff{t}
			\\
			=-\iint_{\omega\times(0,T)} \nabla \left(\zeta e^{-2s\alpha}s^{2^{(m+2)}-1}\xi^{2^{(m+2)}-1}\right)\cdot \nabla \varphi_1^{(1)} \varphi_1^{(1)} \diff{x}\diff{t}
			\\
			-\iint_{\omega\times(0,T)}  \zeta e^{-2s\alpha}s^{2^{(m+2)}-1}\xi^{2^{(m+2)}-1} \Delta \varphi_1^{(1)} \varphi_1^{(1)} \diff{x}\diff{t}.
		\end{multline*}
		
		Considering \eqref{11:01}, and using Young's inequality, we deduce the existence of $C$ such that for any $s\geq C$ and for any $\varepsilon >0,$
		\begin{multline} \label{22:41}
			\iint_{\widetilde\omega\times(0,T)} e^{-2s\alpha}s^{2^{(m+2)}-1}\xi^{2^{(m+2)}-1} \left|\nabla \varphi_1^{(1)}\right|^2 \diff{x}\diff{t}
			\leq 
			\varepsilon
			\iint_Q e^{-2s\alpha}\left( s^6\xi^{6}\left|\nabla\varphi_1^{(1)}\right|^2 
			+ s^5\xi^5\left|\Delta\varphi_1^{(1)}\right|^2 \right)  \diff{x} \diff{t}
			\\
			+\frac{C}{\varepsilon}\iint_{\omega\times(0,T)}e^{-2s\alpha}s^{2^{(m+3)}-6}\xi^{2^{(m+3)}-6} \left| \varphi_1^{(1)} \right|^2  \diff{x} \diff{t}.
		\end{multline}		
Gathering \eqref{eq:finiteiteration}, \eqref{eq:deltavarphi1nablavarphi1} and the above estimate implies that for any $\varepsilon>0$, there exists a constant $C>0$ such that
\begin{multline}\label{8415:56}
\sum_{i=1}^m I(s,\varphi^{(i)})
\leq  \varepsilon \iint_Q e^{-2s\alpha} \left(s^6\xi^{6}\left|\nabla\varphi_1^{(1)}\right|^2 +s^5\xi^5 \left| \Delta\varphi_1^{(1)}\right|^2+s^3\xi^3 \left|\nabla \Delta\varphi_1^{(1)}\right|^2\right) \diff{x}\diff{t}
			\\
			+\frac{C}{\varepsilon}\iint_{\omega\times(0,T)}e^{-2s\alpha}s^{2^{(m+3)}-6}\xi^{2^{(m+3)}-6} \left| \varphi_1^{(1)} \right|^2  \diff{x} \diff{t}.
\end{multline}
On the other hand, applying \cref{lemma:215}, we deduce that
\begin{multline*}
\iint_Q  e^{-2s\alpha}s^8\xi^8 \left|\varphi_1^{(1)}\right|^2 \diff{x}\diff{t}  
			+\iint_Q  e^{-2s\alpha}s^6\xi^6 \left|\nabla \varphi_1^{(1)}\right|^2 \diff{x}\diff{t}  
\\
\leq
C \iint_Q  e^{-2s\alpha}s^5\xi^5 \left|\Delta\varphi_1^{(1)}\right|^2 \diff{x}\diff{t}  
			+\iint_{\omega\times (0,T)}  e^{-2s\alpha}s^8\xi^8 \left|\varphi_1^{(1)}\right|^2 \diff{x}\diff{t}.
\end{multline*}
Combining the above estimate with \eqref{8415:56} and with the definition \eqref{defI1} of $I$, we deduce that
$$
\sum_{i=1}^m I(s,\varphi^{(i)})
\leq \varepsilon I(s,\varphi^{(1)}) 
+ \frac{C}{\varepsilon}\iint_{\omega\times(0,T)}e^{-2s\alpha}s^{2^{(m+3)}-6}\xi^{2^{(m+3)}-6} \left| \varphi_1^{(1)} \right|^2  \diff{x} \diff{t}
$$
and this yields the conclusion of \cref{prop:0}.
\end{proof}

\section{Proof of the main results}\label{S4}
\subsection{Final state observability}
In this section, we use \cref{prop:0} in order to prove the final state observability of the adjoint system \eqref{eq:systemadjoint2}.
\begin{lema}\label{L01}
Assume $T\in (0,1)$ and $\omega$ is non empty open set of $\Omega$. Then, there exists $C>0$ and $\ell \geq 14$ such that
for any $\varphi_0\in \mathcal{H}$, the solution $\varphi$ of \eqref{eq:systemadjoint2} satisfies
    \begin{equation}\label{tt44}
        \sum_{i=1}^m \int_{\Omega}  \left|\varphi^{(i)}(x,0)\right|^2  \diff{x}
        \leq C e^{\frac{C}{T^{2\ell}}}
        \iint_{\omega\times(0,T)} \left| \varphi_1^{(1)} \right|^2  \diff{x} \diff{t}.
    \end{equation}
\end{lema}
\begin{proof}
    First, we consider an energy estimate of the adjoint system \eqref{eq:systemadjoint2}. Multiplying each equation \eqref{tt00} by $\varphi^{(i)}$ and integrating by parts, we deduce 
        \begin{equation*}
    	    -\frac12 \frac{\diff{}}{\diff{t}} \sum_{i=1}^m \int_{\Omega}  \left|\varphi^{(i)}\right|^2  \diff{x}
        	+\sum_{i=1}^m \int_{\Omega} \left|\nabla \varphi^{(i)}\right|^2  \diff{x}
    	    =\sum_{i,j=1}^m \int_{\Omega} \left( A_{j,i}\varphi^{(i)}\cdot \varphi^{(j)} +\left[\left(B_{j,i}\cdot\nabla\right)\varphi^{(i)}\right]\cdot \varphi^{(j)} \right) \diff{x}.
        \end{equation*}
    Thus, using the Gr\"onwall lemma, there exists $C>0$ such that
    $$
    t\mapsto e^{Ct}\sum_{i=1}^m\int_{\Omega}  \left|\varphi^{(i)}(x,t)\right|^2  \diff{x}
    $$
    is nondecreasing. In particular, for some constant $C>0$,
        \begin{equation} \label{tt40}
            \sum_{i=1}^m \int_{\Omega}  \left|\varphi^{(i)}(x,0)\right|^2  \diff{x}
            \leq \frac{2}{T} e^{CT}
            \sum_{i=1}^m \int_{T/4}^{3T/4}  \int_{\Omega}  \left|\varphi^{(i)}(x,t)\right|^2  \diff{x}  \diff{t}.
        \end{equation}

    On the other hand, from \eqref{eq:phiiphi11} and \eqref{defI1}, we deduce that
        \begin{equation}\label{tt41}
            \sum_{i=1}^m \iint_Q e^{-2s\alpha^*}(s\xi^*)^5  \left|\varphi^{(i)} \right|^2 \diff{x}\diff{t}
    	    \leq C\iint_{\omega\times(0,T)}e^{-2s\alpha}\left(s\xi\right)^{2^{m+3}-6} \left| \varphi_1^{(1)} \right|^2  \diff{x} \diff{t}.
        \end{equation}
    Using that for $t\in [T/4,3T/4]$,
    $$
    \frac{3T^2}{16}\leq t(T-t)\leq \frac{T^2}{4},
    $$
    we deduce, from  \eqref{eq:weights2}, the existence of two constants $C_1,C_2>0$ such that
    for $t\in [T/4,3T/4]$,
    $$
    \alpha^*(t) \leq \frac{C_1}{T^{2\ell}}, \quad \xi^*(t) \geq \frac{C_2}{T^{2\ell}}
    $$
    and consequently, for some constant $C_3>0$,
        \begin{equation} \label{eq:cost4}
            e^{-2s\alpha^*}(s\xi^*)^5\geq e^{-\frac{C_3}{T^{2\ell}} }.
        \end{equation}	
    Similarly, from  \eqref{eq:weights}, there exist two constants $c_1,c_2>0$ such that
    for $(x,t)\in \Omega\times [0,T]$,
    $$
    \alpha(x,t) \geq \frac{c_1}{t^\ell(T-t)^\ell}, \quad \xi(x,t) \leq \frac{c_2}{t^\ell(T-t)^\ell}
    $$
    and consequently, for some constant $c_3>0$,
        \begin{equation} \label{eq:cost40}
            e^{-2s\alpha}\left(s\xi\right)^{2^{m+3}-6}
            \leq e^{-2s\frac{c_1}{t^\ell(T-t)^\ell}}\left(s \frac{c_2}{t^\ell(T-t)^\ell}\right)^{2^{m+3}-6}
            \leq c_3.
        \end{equation}	
    Combining \eqref{tt40}, \eqref{tt41}, \eqref{eq:cost4} and \eqref{eq:cost40}, we deduce that for some constant $C$,
    $$
    \sum_{i=1}^m \int_{\Omega} \left|\varphi^{(i)}(x,0)\right|^2  \diff{x}
    \leq 
    \frac{C}{T} e^{C\left(T+\frac{1}{T^{2\ell}}\right) }
    \iint_{\omega\times(0,T)} \left| \varphi_1^{(1)} \right|^2  \diff{x} \diff{t}.
    $$
    This implies that \eqref{tt44}.
\end{proof}

\subsection{Proof of \cref{teo:1}}
We use the functional framework introduced in \cref{sec_reg}.
We recall that $H$ and $V$ are defined by \eqref{defH} and \eqref{defV}, that $P_0 : L^2(\Omega)^N \to H$ is the Leray projector.
We define the control operator $\mathcal{B}\in \mathcal{L}(\mathcal{U},\mathcal{H})$ by
    \begin{equation}\label{defcalUB}
        \mathcal{U}:=L^2(\omega)^{N-1},
        \quad
        \mathcal{B} v = \mathcal{B} 
        \begin{pmatrix}
        v_1,\ldots,v_{N-1}
        \end{pmatrix}
        :=
        \begin{pmatrix}
        P_0\left(\begin{pmatrix}
        v_1,\ldots,v_{N-1},0
        \end{pmatrix}1_{\omega}\right),
        0, \ldots,0
        \end{pmatrix}.
    \end{equation}
With the above definition and the definition \eqref{defcalAdom}-\eqref{defcalA} of $\mathcal{A}$, we can write \eqref{eq:sistemamatriz} as
    \begin{equation} \label{tt50}
	    \begin{cases}
		    y'+\mathcal{A}y=\mathcal{B} v, 
		    \\
		    y(0)=y_0.
	    \end{cases}
    \end{equation}
As it is well-known (see, for instance, \cite[p.357]{TucsnakWeiss}), system \eqref{tt50} is null-controllable in time $T>0$ if and only if 
there exists $K(T)>0$ such that
    \begin{equation}\label{tt51}
        \left\| e^{-T\mathcal{A}^*} \varphi_0\right\|_{\mathcal{H}}^2 \leq K(T)^2 \int_0^T \left\| \mathcal{B}^* e^{-t\mathcal{A}^*} \varphi_0\right\|_{\mathcal{H}}^2  \diff{t}
        \quad (\varphi_0\in \mathcal{H}).
    \end{equation}
Since $\mathcal{A}^*$ is given by \eqref{defcalAstar} and since
$$
\mathcal{B}^* \varphi = \begin{pmatrix}
\left(\varphi^{(1)}_1\right)_{|\omega},\ldots,\left(\varphi^{(1)}_{N-1}\right)_{|\omega}
\end{pmatrix},
$$
we deduce \cref{teo:1} from \cref{L01}.

\subsection{Proof of \cref{teo:2}}
In order to prove \cref{teo:2}, we recall a method introduced in \cite{liu2013single,takahashi2017boundary} to deal with the controllability of nonlinear parabolic systems. We consider $\mathcal{H}$ and $\mathcal{U}$ two Hilbert spaces, $-\mathcal{A} : \mathcal{D}(\mathcal{A}) \to \mathcal{H}$ the infinitesimal generator of an analytic semigroup $\left(e^{-t \mathcal{A}}\right)_{t\geq 0}$ and $B\in \mathcal{L}(\mathcal{U},\mathcal{H})$ a control operator such that 
\eqref{tt51} holds
with $K : (0,\infty)\to [0,\infty)$ continuous and non increasing. 
Let us consider $T>0$ and suppose there exist $\rho_0,\rho_1,\rho \in C^0([0,T],\mathbb{R}^+)$, non increasing, positive in $[0,T)$ such that 
$\rho_0(T)=\rho_1(T)=\rho(T)=0$ and such that, for some constant $q>1$,
    \begin{equation} \label{eq:whoisrho0}
	    \rho_0(t):=\rho_1(q^2(t-T)+T)K((q-1)(T-t))\quad\left(t\in \left[T\left(1-\frac{1}{q^2}\right),T\right] \right),
    \end{equation}
    \begin{equation} \label{eq:whoisrho}
	    \rho_0\leq C\rho,\quad\rho_1\leq C\rho,\quad |\rho'|\rho_0\leq C\rho^2 \quad \left(t\in [0,T]\right).
    \end{equation}
for some constant $C>0$.
We denote by $L^2_{\rho_1}(0,T;\mathcal{H})$ the space
$$
L^2_{\rho_1}(0,T;\mathcal{H}):=\left\{f\in L^2(0,T;\mathcal{H}) \ ; \ \frac{f}{\rho_1}\in L^2(0,T;\mathcal{H})\right\}
$$
and we define similarly $L^2_{\rho_0}(0,T;\mathcal{U})$.

Then we can consider the control problem
\begin{equation} \label{tt500}
	\begin{cases}
		y'+\mathcal{A}y=\mathcal{B} v+f, 
		\\
		y(0)=y_0.
	\end{cases}
\end{equation}
We have the following result (see \cite{liu2013single}):
\begin{teorema}\label{Tliu}
    With the above assumptions, there exists a bounded operator 
    $$
    \mathcal{E}_T\in \mathcal{L}\left( \mathcal{D}\left(\mathcal{A}^{1/2}\right)\times L^2_{\rho_1}(0,T;\mathcal{H}), L^2_{\rho_0}(0,T;\mathcal{U}) \right)
    $$
    such that for any $y_0\in \mathcal{D}\left(\mathcal{A}^{1/2}\right)$
    and for any $f\in L^2_{\rho_1}(0,T;\mathcal{H})$,  
    the solution $y$ of \eqref{tt500} with $u=\mathcal{E}_T(y_0,f)$ satisfies
    $$
    \frac{y}{\rho} \in L^2(0,T;\mathcal{D}(\mathcal{A}))\cap C^0\left([0,T];\mathcal{D}\left(\mathcal{A}^{1/2}\right)\right) \cap H^1(0,T;\mathcal{H}).
    $$
    Moreover there exists a constant $C$ such that 
    $$
    \left\| \frac{y}{\rho} \right\|_{L^2(0,T;\mathcal{D}(\mathcal{A}))\cap C^0\left([0,T];\mathcal{D}\left(\mathcal{A}^{1/2}\right)\right) \cap H^1(0,T;\mathcal{H})}
    \leq C\left( \|y_0\|_{\mathcal{D}\left(\mathcal{A}^{1/2}\right)}+\left\| f\right\|_{L^2_{\rho_1}(0,T;\mathcal{H})} \right).
    $$
\end{teorema}

\begin{observacion}
Note that in \cite{liu2013single}, $\mathcal{A}$ is assumed to be self-adjoint positive but the result can be extended to the case where $-\mathcal{A}$ is the generator of an analytic semigroup. Indeed, the hypothesis used in the proof is the maximal regularity of \eqref{tt500} for $v=0$.
\end{observacion}

\begin{observacion}
Since $\rho(T)=0$, the above result implies in particular that $y(T)=0$, that is the null-controllability of \eqref{tt500}.
\end{observacion}

In the previous section, we have defined for our problem the spaces $\mathcal{H}$, $\mathcal{U}$, $\mathcal{A}$ and $\mathcal{B}$,
see \eqref{defcalH}, \eqref{defcalA} and \eqref{defcalUB}. We have shown in \cref{sec_reg} that $-\mathcal{A}$ is the generator of an analytic semigroup.
%
Finally, applying \cref{L01}, we deduce that \eqref{tt51} holds with 
$$
K(T)=C_K e^{\frac{C_K}{T^{2\ell}}},
$$
for some constant $C_K>0$. 

Let us consider 
$$
q \in \left(1,2^{\frac{1}{4\ell}}\right).
$$
and let us set
$$
\rho_0(t):=C_K e^{-\frac{C_0}{(T-t)^{2\ell}}}, \quad \rho_1(t):=e^{-\frac{C_1}{(T-t)^{2\ell}}}, \quad
\rho(t):=e^{-\frac{C_\star}{(T-t)^{2\ell}}}
$$ 
with $C_0,C_1,C_\star$ some positive constants such that
$$
C_0:=\frac{C_1}{q^{4\ell}}-\frac{C_K}{(q-1)^{2\ell}} > \frac{C_1}{2}, \quad 
\frac{C_1}{2}<C_\star<C_0<C_1.
$$
Then we can check that \eqref{eq:whoisrho0} and \eqref{eq:whoisrho} hold and we have moreover that
\begin{equation}\label{tt52}
    \rho^2 \leq \rho_1.
\end{equation}
Consequently, we deduce from \cref{Tliu} a controllability result on the system
    \begin{equation} \label{eq:system+f}
    \left\{\begin{array}{lll}
    	\partial_t y^{(1)}-\Delta y^{(1)}+\nabla p^{(1)}=\sum_{j=1}^{m}\left( B_{1,j}\cdot\nabla\right) y^{(j)}+\sum_{j=1}^{m} A_{1,j} y^{(j)}+v e_1 1_{\omega}+f^{(1)} &\text{in} \ Q, & \\
	    \partial_t y^{(i)}-\Delta y^{(i)}+\nabla p^{(i)}=\sum_{j=i}^{m}\left( B_{i,j}\cdot\nabla\right) y^{(j)}+\sum_{j=i-1}^m A_{i,j} y^{(j)} +f^{(i)} &\text{in} \ Q, &(2\leq i\leq m)\\
    	\nabla\cdot y^{(i)}=0 &\text{in }Q, &(1\leq i\leq m)\\
	    y^{(i)}=0 &\text{on }\Sigma, &(1\leq i\leq m)\\
    	y^{(i)}(\cdot,0)=y_0^{(i)} &\text{in }\Omega. &(1\leq i\leq m)
    \end{array}\right.
    \end{equation}	
More precisely, there exists
$$
\mathcal{E}_T\in \mathcal{L}\left( \mathcal{V} \times L^2_{\rho_1}(0,T;\left[L^2(\Omega)^{N}\right]^m), L^2_{\rho_0}(0,T;L^2(\omega)) \right)
$$
such that for any $y_0\in \mathcal{V}$ and for any $f=\begin{pmatrix}
f^{(1)},\ldots,f^{(m)}
\end{pmatrix}\in L^2_{\rho_1}(0,T;\left[L^2(\Omega)^{N}\right]^m)$, the solution $y$ of \eqref{eq:system+f} 
with the control $v=\mathcal{E}_T(y_0,f)$ satisfies
    \begin{equation}\label{tt60}
        \frac{y}{\rho} \in L^2(0,T;\left[H^2(\Omega)^N\right]^m)\cap C^0\left([0,T];\left[H^1(\Omega)^N\right]^m\right) \cap H^1(0,T;\mathcal{H}).
    \end{equation}
Moreover we have the following estimate 
    \begin{equation}\label{tt61}
        \left\| \frac{y}{\rho} \right\|_{L^2(0,T;\left[H^2(\Omega)^N\right]^m)\cap C^0\left([0,T];\left[H^1(\Omega)^N\right]^m\right) \cap H^1(0,T;\mathcal{H})}
        \leq C\left( \|y_0\|_{\mathcal{V}}+\left\| f\right\|_{L^2_{\rho_1}(0,T;\left[L^2(\Omega)^N\right]^m)} \right).
    \end{equation}

We are now in a position to prove \cref{teo:2}:
\begin{proof}[Proof of \cref{teo:2}]
    First we notice that $y$ is solution of \eqref{eq:sistemamatriz+navier} if it is a solution of \eqref{eq:system+f} with
    $$
    f=\left(-(y^{(1)}\cdot\nabla)y^{(1)},...,-(y^{(m)}\cdot\nabla)y^{(m)}\right).
    $$
    Thus, we consider the mapping
    $$
    \mathcal{N}_T: f\in B_R \mapsto (-(y^{(1)}\cdot\nabla)y^{(1)},...,-(y^{(m)}\cdot\nabla)y^{(m)}),
    $$
    where
    $$
    B_R:=\left\{ f\in L^2_{\rho_1}(0,T;\left[L^2(\Omega)^N\right]^m) \ ; \ \left\|\frac{f}{\rho_1}\right\|_{L^2(0,T;\left[L^2(\Omega)^N\right]^m)} \leq R\right\}
    $$
    where $R>0$ is such that
    $$
    \|y_0\|_{\mathcal{V}}\leq R.
    $$
    We are going to show that for $R$ small enough (and thus $\|y_0\|_{\mathcal{V}}$ small enough), $\mathcal{N}_T(B_R)\subset B_R$ and that 
    $\left(\mathcal{N}_T\right)_{|B_R}$ is a strict contraction. Using the Banach fixed point theorem we deduce the existence of a fixed point of $\mathcal{N}_T$.
    The corresponding solution $y$ of \eqref{eq:system+f} is a solution of \eqref{eq:sistemamatriz+navier} and from \eqref{tt60}, we deduce that
    $y(\cdot,T)=0$.

    It thus remains to prove that for $R$ small enough, $\mathcal{N}_T(B_R)\subset B_R$ and that 
    $\left(\mathcal{N}_T\right)_{|B_R}$ is a strict contraction. 
    In order to do this, we first note that, using \eqref{tt52}, Sobolev's embeddings and H\"older's inequalities, we have
        \begin{multline}\label{tt62}
            \int_0^T\int_\Omega \left |\frac{v\cdot\nabla w}{\rho_1}\right|^2 \diff{x} \diff{t}
            \leq 
            \int_0^T\int_\Omega \left | \left(\frac{v}{\rho}\right)\cdot\nabla \left(\frac{w}{\rho}\right)\right|^2 \diff{x} \diff{t}
            \leq C\bigg\|\frac{v}{\rho}\bigg\|^2_{L^\infty(0,T;L^6(\Omega)^N)}\bigg\|\frac{\nabla w}{\rho}\bigg\|^2_{L^2(0,T;L^6(\Omega)^N)} 
            \\
            \leq C\bigg\|\frac{v}{\rho}\bigg\|^2_{L^\infty(0,T;H^1(\Omega)^N)} \bigg\|\frac {w}{\rho}\bigg\|^2_{L^2(0,T;H^2(\Omega)^N)}.
        \end{multline}
    Using this relation and \eqref{tt61}, we deduce that
    $$
    \left\| \frac{\mathcal{N}_T(f)}{\rho_1} \right\|_{L^2(0,T;\left[L^2(\Omega)^N\right]^m)} \leq 
    C\left( \|y_0\|_{\mathcal{V}}+\left\| f\right\|_{L^2_{\rho_1}(0,T;\left[L^2(\Omega)^N\right]^m)} \right)^2
    \leq 4CR^2\leq R,
    $$
    for $R$ small enough. For such $R$, we have $\mathcal{N}_T(B_R)\subset B_R$.

    Now, let us consider $\widetilde f, \widehat f\in B_R$ and let us write $f=\widetilde f- \widehat f$.
    We consider the solution $\widetilde y$ (resp. $\widehat y$) the solution of \eqref{eq:system+f} 
    associated with the control $\widetilde v=\mathcal{E}_T(y_0, \widetilde f)$ (resp. $\widehat v=\mathcal{E}_T(y_0, \widehat f)$).
    Then, $y:=\widetilde y- \widehat y$ is the solution of \eqref{eq:system+f} 
    associated with the control $v:=\mathcal{E}_T(0,f)$ and thus
    $$
    \left\| \frac{y}{\rho} \right\|_{L^2(0,T;\left[H^2(\Omega)^N\right]^m)\cap C^0\left([0,T];\left[H^1(\Omega)^N\right]^m\right) \cap H^1(0,T;\mathcal{H})}
    \leq C\left\| f\right\|_{L^2_{\rho_1}(0,T;\left[L^2(\Omega)^N\right]^m)}.
    $$
    Using this and \eqref{tt62}, we obtain
        \begin{multline*}
        \left\| \frac{\mathcal{N}_T(\widetilde f)}{\rho_1}-\frac{\mathcal{N}_T(\widehat f)}{\rho_1} \right\|_{L^2(0,T;\left[L^2(\Omega)^N\right]^m)}
        \leq
        \left\| \left(\frac{\widetilde y}{\rho}\right)\cdot\nabla \left(\frac{y}{\rho}\right) \right\|_{L^2(0,T;\left[L^2(\Omega)^N\right]^m)}+
        \left\| \left(\frac{y}{\rho}\right)\cdot\nabla \left(\frac{\widehat y}{\rho}\right) \right\|_{L^2(0,T;\left[L^2(\Omega)^N\right]^m)}
        \\
        \leq CR \left\| f\right\|_{L^2_{\rho_1}(0,T;\left[L^2(\Omega)^N\right]^m)}.
        \end{multline*}
        Thus for $R$ small enough, $\left(\mathcal{N}_T\right)_{|B_R}$ is a strict contraction and this ends the proof of \cref{teo:2}.
    \end{proof}

\bibliographystyle{abbrv}
\bibliography{Bibliografia}

\begin{thebibliography}{10}

\bibitem{Ammarkhodja2011}
F.~Ammar-Khodja, A.~Benabdallah, M.~Gonz\'{a}lez-Burgos, and L.~de~Teresa.
\newblock Recent results on the controllability of linear coupled parabolic
  problems: a survey.
\newblock {\em Math. Control Relat. Fields}, 1(3):267--306, 2011.

\bibitem{Ammarkhodja2014}
F.~Ammar~Khodja, A.~Benabdallah, M.~Gonz\'{a}lez-Burgos, and L.~de~Teresa.
\newblock Minimal time for the null controllability of parabolic systems: the
  effect of the condensation index of complex sequences.
\newblock {\em J. Funct. Anal.}, 267(7):2077--2151, 2014.

\bibitem{Ammarkhodja2016}
F.~Ammar~Khodja, A.~Benabdallah, M.~Gonz\'{a}lez-Burgos, and L.~de~Teresa.
\newblock New phenomena for the null controllability of parabolic systems:
  minimal time and geometrical dependence.
\newblock {\em J. Math. Anal. Appl.}, 444(2):1071--1113, 2016.

\bibitem{Benabdallah2020}
A.~Benabdallah, F.~Boyer, and M.~Morancey.
\newblock A block moment method to handle spectral condensation phenomenon in
  parabolic control problems.
\newblock {\em Ann. H. Lebesgue}, 3:717--793, 2020.

\bibitem{MR3348416}
N.~Carre\~{n}o, S.~Guerrero, and M.~Gueye.
\newblock Insensitizing controls with two vanishing components for the
  three-dimensional {B}oussinesq system.
\newblock {\em ESAIM Control Optim. Calc. Var.}, 21(1):73--100, 2015.

\bibitem{Carreno-Gueye2014}
N.~Carre\~{n}o and M.~Gueye.
\newblock Insensitizing controls with one vanishing component for the
  {N}avier-{S}tokes system.
\newblock {\em J. Math. Pures Appl. (9)}, 101(1):27--53, 2014.

\bibitem{Conforto2017}
F.~Conforto, L.~Desvillettes, and R.~Monaco.
\newblock Some asymptotic limits of reaction-diffusion systems appearing in
  reversible chemistry.
\newblock {\em Ric. Mat.}, 66(1):99--111, 2017.

\bibitem{MR1470445}
J.-M. Coron and A.~V. Fursikov.
\newblock Global exact controllability of the {$2$}{D} {N}avier-{S}tokes
  equations on a manifold without boundary.
\newblock {\em Russian J. Math. Phys.}, 4(4):429--448, 1996.

\bibitem{coron2009null}
J.-M. Coron and S.~Guerrero.
\newblock Null controllability of the {$N$}-dimensional {S}tokes system with
  {$N-1$} scalar controls.
\newblock {\em J. Differential Equations}, 246(7):2908--2921, 2009.

\bibitem{MR3279537}
J.-M. Coron and P.~Lissy.
\newblock Local null controllability of the three-dimensional {N}avier-{S}tokes
  system with a distributed control having two vanishing components.
\newblock {\em Invent. Math.}, 198(3):833--880, 2014.

\bibitem{ErdiToth1989}
P.~\'{E}rdi and J.~T\'{o}th.
\newblock {\em Mathematical models of chemical reactions}.
\newblock Nonlinear Science: Theory and Applications. Princeton University
  Press, Princeton, NJ, 1989.
\newblock Theory and applications of deterministic and stochastic models.

\bibitem{fattorini1975}
H.~O. Fattorini and D.~L. Russell.
\newblock Uniform bounds on biorthogonal functions for real exponentials with
  an application to the control theory of parabolic equations.
\newblock {\em Quart. Appl. Math.}, 32:45--69, 1974/75.

\bibitem{MR2224822}
E.~Fern\'{a}ndez-Cara, M.~Gonz\'{a}lez-Burgos, S.~Guerrero, and J.-P. Puel.
\newblock Null controllability of the heat equation with boundary {F}ourier
  conditions: the linear case.
\newblock {\em ESAIM Control Optim. Calc. Var.}, 12(3):442--465, 2006.

\bibitem{fernandez2004local}
E.~Fern\'{a}ndez-Cara, S.~Guerrero, O.~Y. Imanuvilov, and J.-P. Puel.
\newblock Local exact controllability of the {N}avier-{S}tokes system.
\newblock {\em J. Math. Pures Appl. (9)}, 83(12):1501--1542, 2004.

\bibitem{fernandez2006some}
E.~Fern\'{a}ndez-Cara, S.~Guerrero, O.~Y. Imanuvilov, and J.-P. Puel.
\newblock Some controllability results for the {$N$}-dimensional
  {N}avier-{S}tokes and {B}oussinesq systems with {$N-1$} scalar controls.
\newblock {\em SIAM J. Control Optim.}, 45(1):146--173, 2006.

\bibitem{fursikov1996controllability}
A.~V. Fursikov and O.~Y. Imanuvilov.
\newblock {\em Controllability of evolution equations}, volume~34 of {\em
  Lecture Notes Series}.
\newblock Seoul National University, Research Institute of Mathematics, Global
  Analysis Research Center, Seoul, 1996.

\bibitem{gonzalez2010controllability}
M.~Gonz\'{a}lez-Burgos and L.~de~Teresa.
\newblock Controllability results for cascade systems of {$m$} coupled
  parabolic {PDE}s by one control force.
\newblock {\em Port. Math.}, 67(1):91--113, 2010.

\bibitem{SergioNavier}
S.~Guerrero.
\newblock Local exact controllability to the trajectories of the
  {Navier}-{Stokes} system with nonlinear {Navier}-slip boundary conditions.
\newblock {\em ESAIM, Control Optim. Calc. Var.}, 12:484--544, 2006.

\bibitem{Guerrero2007}
S.~Guerrero.
\newblock Controllability of systems of {S}tokes equations with one control
  force: existence of insensitizing controls.
\newblock {\em Ann. Inst. H. Poincar\'{e} Anal. Non Lin\'{e}aire},
  24(6):1029--1054, 2007.

\bibitem{guerrero2007null}
S.~Guerrero.
\newblock Null controllability of some systems of two parabolic equations with
  one control force.
\newblock {\em SIAM J. Control Optim.}, 46(2):379--394, 2007.

\bibitem{SergioCristhian}
S.~Guerrero and C.~Montoya.
\newblock Local null controllability of the {{\(N\)}}-dimensional
  {N}avier-{S}tokes system with nonlinear {N}avier-slip boundary conditions and
  {{\(N-1\)}} scalar controls.
\newblock {\em J. Math. Pures Appl. (9)}, 113:37--69, 2018.

\bibitem{Iida2017}
M.~Iida, H.~Monobe, H.~Murakawa, and H.~Ninomiya.
\newblock Vanishing, moving and immovable interfaces in fast reaction limits.
\newblock {\em J. Differential Equations}, 263(5):2715--2735, 2017.

\bibitem{imanuvilov1998}
O.~Y. Imanuvilov.
\newblock On exact controllability for the {N}avier-{S}tokes equations.
\newblock {\em ESAIM Control Optim. Calc. Var.}, 3:97--131, 1998.

\bibitem{imanuvilov2001remarks}
O.~Y. Imanuvilov.
\newblock Remarks on exact controllability for the {N}avier-{S}tokes equations.
\newblock {\em ESAIM Control Optim. Calc. Var.}, 6:39--72, 2001.

\bibitem{imanuvilov2009carleman}
O.~Y. Imanuvilov, J.-P. Puel, and M.~Yamamoto.
\newblock Carleman estimates for parabolic equations with nonhomogeneous
  boundary conditions.
\newblock {\em Chin. Ann. Math. Ser. B}, 30(4):333--378, 2009.

\bibitem{lebeau1995}
G.~Lebeau and L.~Robbiano.
\newblock Contr\^{o}le exacte de l'\'{e}quation de la chaleur.
\newblock In {\em S\'{e}minaire sur les \'{E}quations aux {D}\'{e}riv\'{e}es
  {P}artielles, 1994--1995}, pages Exp. No. VII, 13. \'{E}cole Polytech.,
  Palaiseau, 1995.

\bibitem{lionsnon1}
J.-L. Lions and E.~Magenes.
\newblock {\em Non-homogeneous boundary value problems and applications. {V}ol.
  {I}}.
\newblock Die Grundlehren der mathematischen Wissenschaften, Band 181.
  Springer-Verlag, New York-Heidelberg, 1972.
\newblock Translated from the French by P. Kenneth.

\bibitem{lionsnon2}
J.-L. Lions and E.~Magenes.
\newblock {\em Non-homogeneous boundary value problems and applications. {V}ol.
  {II}}.
\newblock Die Grundlehren der mathematischen Wissenschaften, Band 182.
  Springer-Verlag, New York-Heidelberg, 1972.
\newblock Translated from the French by P. Kenneth.

\bibitem{liu2013single}
Y.~Liu, T.~Takahashi, and M.~Tucsnak.
\newblock Single input controllability of a simplified fluid-structure
  interaction model.
\newblock {\em ESAIM Control Optim. Calc. Var.}, 19(1):20--42, 2013.

\bibitem{MontoyadeTeresa2018}
C.~Montoya and L.~de~Teresa.
\newblock Robust {S}tackelberg controllability for the {N}avier-{S}tokes
  equations.
\newblock {\em NoDEA Nonlinear Differential Equations Appl.}, 25(5):Paper No.
  46, 33, 2018.

\bibitem{Pazy}
A.~Pazy.
\newblock {\em Semigroups of linear operators and applications to partial
  differential equations}, volume~44 of {\em Applied Mathematical Sciences}.
\newblock Springer-Verlag, New York, 1983.

\bibitem{Sohr}
H.~Sohr.
\newblock {\em The {N}avier-{S}tokes equations}.
\newblock Modern Birkh\"{a}user Classics. Birkh\"{a}user/Springer Basel AG,
  Basel, 2001.
\newblock An elementary functional analytic approach, [2013 reprint of the 2001
  original] [MR1928881].

\bibitem{takahashi2017boundary}
T.~Takahashi.
\newblock Boundary local null-controllability of the {K}uramoto-{S}ivashinsky
  equation.
\newblock {\em Math. Control Signals Systems}, 29(1):Art. 2, 21, 2017.

\bibitem{Temam}
R.~Temam.
\newblock {\em Navier-{S}tokes equations}, volume~2 of {\em Studies in
  Mathematics and its Applications}.
\newblock North-Holland Publishing Co., Amsterdam-New York, revised edition,
  1979.
\newblock Theory and numerical analysis, With an appendix by F. Thomasset.

\bibitem{TucsnakWeiss}
M.~Tucsnak and G.~Weiss.
\newblock {\em Observation and control for operator semigroups}.
\newblock Birkh\"{a}user Advanced Texts: Basler Lehrb\"{u}cher. [Birkh\"{a}user
  Advanced Texts: Basel Textbooks]. Birkh\"{a}user Verlag, Basel, 2009.

\end{thebibliography}
	
\end{document}